\documentclass{article}

\usepackage{arxiv}

\usepackage[utf8]{inputenc} 
\usepackage[T1]{fontenc}    
\usepackage{hyperref}       
\usepackage{url}            
\usepackage{booktabs}       
\usepackage{amsfonts}       
\usepackage{nicefrac}       
\usepackage{microtype}      
\usepackage{doi}
\usepackage{style}
 \usepackage{enumitem}
 
\title{Two continuous extensions of the Neural Approximated Virtual Element Method}

\author{ \href{https://orcid.org/0000-0001-8642-4258}{\includegraphics[scale=0.06]{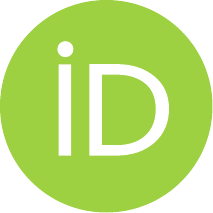}\hspace{1mm}Stefano~Berrone} \\
	Dipartimento di Scienze Matematiche\\
	``G. L. Lagrange''\\
	Politecnico di Torino, TO, 10129 \\
	\texttt{stefano.berrone@polito.it} \\
	\And
	\href{https://orcid.org/0000-0002-2837-2792}{\includegraphics[scale=0.06]{orcid.pdf}\hspace{1mm}Moreno~Pintore} \\
	Laboratoire Jacques-Louis Lions, \\
	Sorbonne Universit\'e, \\
	MEGAVOLT Team, Inria, \\
	4 place Jussieu, 75005 Paris, France \\
	\texttt{moreno.pintore@sorbonne-universite.fr} \\
	\And
	\href{https://orcid.org/0000-0002-8540-3639}{\includegraphics[scale=0.06]{orcid.pdf}\hspace{1mm}Gioana~Teora} \\
	Dipartimento di Scienze Matematiche\\
	``G. L. Lagrange''\\
	Politecnico di Torino, TO, 10129 \\
	\texttt{gioana.teora@polito.it} \\
}

\renewcommand{\shorttitle}{Two continuous extensions of the Neural Approximated Virtual Element Method}

\hypersetup{
pdftitle={Two continuous extensions of the Neural Approximated Virtual Element Method},
pdfsubject={q-bio.NC, q-bio.QM},
pdfauthor={Stefano~Berrone, Moreno~Pintore, Gioana~Teora},
pdfkeywords={Virtual Element Method, Neural Network, Basis Functions, Polygonal Meshes},
}

\begin{document}
\maketitle

\begin{abstract}
We propose two globally continuous neural-based variants of the Neural Approximated Virtual Element Method (NAVEM), termed B-NAVEM and P-NAVEM. Both approaches construct local basis functions using pre-trained fully connected neural networks while ensuring exact continuity across adjacent mesh elements. B-NAVEM leverages a Physics-Informed Neural Network to approximately solve the local Laplace problem that defines the virtual element basis functions, whereas P-NAVEM directly enforces polynomial reproducibility via a tailored loss function, without requiring harmonicity within the element interior. Numerical experiments assess the methods in terms of computational cost, memory usage, and accuracy during both training and testing phases.
\end{abstract}


\section{Introduction}

The interest in Galerkin methods for the approximation of solutions to Partial Differential Equations (PDEs) based on polytopal (i.e. polygonal and polyhedral) meshes has recently
grown in the last decade, due to the high flexibility that such meshes allow. Indeed, the usage of polytopal meshes automatically includes the possibility of using non-convex elements, hanging nodes (enabling natural
handling of interface problems and refinement strategies with non-matching grids), easier construction of adaptive meshes, and
efficient approximations of geometric data features \cite{Benedetto2016, Vicini2024, GrappeinTeora2025}.
Examples of polytopal methods include Hybrid High-Order methods (HHO) \cite{inertiaHHO}, Hybrid Discontinuous Galerkin methods (HDGM) \cite{Cockburn2009}, and polygonal discontinuous Galerkin methods (DG-FEM) \cite{Fumagalli2024}.

The Virtual Element Method (VEM) is an alternative 
approach that enables discretization on polytopal meshes \cite{LBe16, DaVeigaDassi2017}. It is based on
globally continuous discretization spaces whose trial and test functions are implicitly defined as solutions of a local Laplace problem on each element. For the lowest-order, the elemental VEM space is described as
\begin{align}
     \label{eq:vem:harmonicity}\Vh[E] = \Big\{ v \in \sob{1}{E}:\quad &(i)\ \Delta v = 0, \\ 
     \label{eq:vem:linearity}&(ii) \ v \in \con{0}{\partial E} : v_{|e}\in \Poly{1}{e},\ \forall e \subset \partial E\Big\}, 
\end{align}
where we adopt standard conventions for functional spaces, $E$ and $\partial E$ denote the polygonal element and its boundary, respectively.
Within the virtual element framework, these local PDEs problems are never solved neither exactly nor approximately. As a consequence, the bilinear form, and thereby the entries of the stiffness matrix, are not directly
computable. The computable version relies on an approximate discrete bilinear form
consisting of two additive parts: a polynomial projection-based form ensuring the polynomial patch-test, and a computable stabilizing bilinear form. Despite VEM remarkable flexibility, the usage of stability and projector operators may lead to several issues, especially in the presence of anisotropic problems, eigenvalue problems, nonlinear problems, or when post-processing the quantities of interest \cite{Credali2024, navemElasticity}. For these reasons, new methods have been introduced to approximate VEM functions while retaining VEM advantages and mitigating its limitations. Examples include the Reduced-Basis Virtual Element Method (rbVEM) \cite{Credali2024, Credali2025}, the Lighting Virtual Element Methods (L-VEM) \cite{TrezziZerbinati2024, TrezziZerbinati2025self}, and the Neural Approximated Virtual Element Method (NAVEM) \cite{PintoreTeora2024, PintoreTeora2025}. The last two approaches describe virtual functions as a linear combination of proper harmonic functions: L-VEM determines these coefficients by using the Laplace Solver \cite{Gopal2019}, while NAVEM leverages neural networks to exploit an efficient offline-online splitting strategy.

Since the elemental approximations in both L-VEM and NAVEM are computed independently on each mesh element, the resulting basis functions are not continuous across adjacent elements. Nevertheless, these methods maintain a $\con{0}{}$-conforming framework without doubling the degrees of freedom along element interfaces. Because virtual basis functions are known in closed form on the boundaries of elements, no approximation is required there, which prevents duplication of boundary degrees of freedom.

In this paper, we propose two variants of NAVEM that restore the continuity of basis functions across neighboring elements. The involved neural networks retain the same kind of architecture as the standard NAVEM, namely a fully-connected feed-forward neural network \cite{goodfellow2016deep}, but differ in the definition of their input/output and loss function. 

In particular, the first variant, referred to as B-NAVEM, employs a Physics-Informed Neural Network (PINN) with exact enforcement of boundary Dirichlet conditions \cite{Sukumar2021, Pintore2023} to approximate the local Laplace problem \eqref{eq:vem:harmonicity}-\eqref{eq:vem:linearity} that defines the VEM Lagrange basis functions. In this setting, the loss function aims to minimize the residual associated with the elemental PDE residual, following the standard PINN paradigm. The second variant, called P-NAVEM, adopts the same technique to enforce boundary Dirichlet data as B-NAVEM, but employs a different loss function that aims at reducing the polynomial reproducibility error. The main goal here is to recover the key property required to guarantee optimal polynomial convergence rates. We emphasize that, in both variants, the basis functions are exact at the boundary of mesh elements, unlike NAVEM, but the functions belonging to these new spaces are no longer exactly harmonic functions in the interior of the elements, like VEM and NAVEM functions. More precisely, B-NAVEM basis functions are harmonic up to a precision that depends on the neural network accuracy, whereas in P-NAVEM this property is neither enforced exactly nor approximately. The proposed two variants are compared with the standard NAVEM in terms of memory usage, computational time, and accuracy both in the training and in the testing phase, through a series of numerical experiments. Moreover, different numerical experiments are proposed to evaluate and compare the performance of the two new neural-based methods against the standard NAVEM and VEM method when solving both linear and nonlinear partial differential problems.

The outline of the paper is as follows. Section~\ref{sec:model_problem_and_nn_arch} presents the neural network strategy that is shared by the three approaches to solve a simple Poisson Problem. In Section~\ref{sec:navem}, we detail the main features of the standard NAVEM construction. Section~\ref{sec:boundary_funcs} describes the strategy adopted to modify neural network output in order to exactly enforce the boundary Dirichlet data, which underlies the B-NAVEM strategy, detailed in Section~\ref{sec:bnavem} and the P-NAVEM approach, which is described in Section~\ref{sec:napem}. Section~\ref{sec:numeri_experiements} reports a set of numerical experiments comparing the different NAVEM strategies and highlighting their advantages with respect to the classical VEM discretization. Finally, Section~\ref{sec:conclusion} summarizes the main findings and outlines future research directions.

\section{Neural networks to approximate basis functions}\label{sec:model_problem_and_nn_arch}

In the following, the superscript $\NN$ will refer to a generic neural-based method, whereas the superscripts $\NAVEM,\ \BNAVEM$ and $\NAPEM$ will denote the specific underlying method, namely the standard NAVEM (or H-NAVEM), B-NAVEM, and P-NAVEM, respectively.

Given a polygonal domain $\Omega\in\R^2$ with boundary $\partial \Omega$, let us consider the following Poisson problem 
\begin{equation}
    \begin{cases}
        -\Delta u = f & \text{in }\Omega,\\
        u  = 0&\text{on }\partial\Omega,
    \end{cases}
    \label{eq:continuous_poisson}
\end{equation}
where the source term $f\in L^2(\Omega)$.  

Let $\V = \sob[0]{1}{\Omega}$, the variational formulation of Problem \eqref{eq:continuous_poisson} reads as: \textit{Find $u\in \V$ such that}
\begin{equation}
    a(u, v) = F(v) \quad \forall v\in \V,
    \label{eq:exact_var_pb}
\end{equation}
where the bilinear form $a:\V\times \V\rightarrow\R$ and the linear form $F:\V\rightarrow \R$ are defined as:
\begin{gather*}
    a(v, w) = \int_\Omega \nabla v\cdot\nabla w, \quad \forall v, w\in \V,\\
F(v) = \int_\Omega f \,v, \quad \forall v\in \V.
\end{gather*}

Let us introduce a tesselation $\Th$ of $\Omega$ made up by polygonal elements $E$. The number of vertices of the polygon $E$ is denoted by $\Nv[E]$. 
Since we consider the lowest-order version, the local virtual element degrees of freedom are the values of functions at the vertices of the element, thus the number of local degrees of freedom $\Ndof[E] = \Nv[E]$. Moreover, the global degrees of freedom are related to the vertices of the tesselation $\Th$ that do not belong to the Dirichlet boundary, and their number is 
\begin{equation}
    \Ndof = \dim \Vh{},
\end{equation}
where the global virtual element space $\Vh{}$ is described as
\begin{equation*}
    \Vh{} = \{v \in \V \cap \con{0}{\overline{\Omega}}:\ v_{|E}  \in \Vh[E]{} \forall E \in \Th \},
\end{equation*}
whereas the local virtual element space $\Vh[E]{}$ is defined in \eqref{eq:vem:harmonicity}-\eqref{eq:vem:linearity}. 

As in a standard Galerkin method, the basic idea of a NAVEM method is to define a finite-dimensional subspace $\nVh{} \subset \V$ as the span of some basis functions that are known in a closed form and approximate Lagrange virtual basis functions spanning $\Vh{}$. In particular, we note that by satisfying 
\begin{itemize}
    \item Property \eqref{eq:vem:linearity}, we gain the $\con{0}{}$-conformity, i.e. continuity across adjacent elements of $\Th$;
    \item Properties \eqref{eq:vem:harmonicity} and \eqref{eq:vem:linearity}, we obtain that $\Poly{k}{E} \subset \nVh{}$. This \textit{polynomial reproducibility} property allows VEM to recover optimal polynomial error convergence estimates.
\end{itemize}

In particular, in the following, we approximate the VEM basis functions using neural networks that aim to learn the following nonlinear map
\begin{equation}
    \left(j, E\right) \to \left(\varphi^{\NN}_{j,E}, \qq^{\NN}_{j,E}\right), \quad \forall j =1,\dots,\Nv[E],\quad \forall E \in \Th.
    \label{eq:nn_map}
\end{equation}
where the input 
\begin{itemize}
    \item the input $\left(j, E\right)$ represents a pair made up by the index $j$ of the elemental degree of freedom (physically represented by a vertex of $E$) and the element $E \in \Th$ itself, identifying the elemental Lagrange basis functions we want to approximate.
    \item the output $\left(\varphi^{\NN}_{j,E}, \qq^{\NN}_{j,E}\right)$ represents the approximation of the virtual basis function $\varphi_{j,E}$ related to the input pair $\left(j, E\right)$ and of its gradient $\nabla \varphi_{j,E}$, respectively.
\end{itemize}
Given the local approximations $\left(\varphi^{\NN}_{j,E}, \qq^{\NN}_{j,E}\right)$ for each pair of input, we can define the local NAVEM space as
\begin{equation}
    \nVh[E]{} = \myspan\{\varphi^{\NN}_{j,E}:\ j =1,\dots,\Ndof[E]\}, \quad \text{and}\quad \nabla \nVh[E]{} = \myspan\{\qq^{\NN}_{j,E}:\ j =1,\dots,\Ndof[E]\}.
    \label{eq:local_navem_space}
\end{equation}
We can now use standard gluing techniques to produce a full set of global basis functions and of their gradients, i.e. $\left(\varphi^{\NN}_{i}, \qq^{\NN}_{i}\right)$, for each $i =1,\dots,\Ndof$ and define the NAVEM global spaces as
\begin{equation*}
    \nVh{} = \myspan\{\varphi^{\NN}_{i}:\ i =1,\dots,\Ndof\}, \quad \text{and}\quad \nabla \nVh{} = \myspan\{\qq^{\NN}_{i}:\ i =1,\dots,\Ndof\}.
\end{equation*}
Finally, we can proceed as in a standard finite element method, and solve the following discrete problem: \textit{Find $u_h\in \nVh{}$ such that}
\begin{equation}
    \sum_{i = 1}^{\Ndof} u_i \sum_{E \in \Th} \int_E \qq^{\NN}_i \cdot \qq^{\NN}_j = \sum_{E \in \Th} \int_E f \varphi^{\NN}_j, \quad \forall j = 1,\dots,\Ndof.
    \label{eq:discrete_var_prob}
\end{equation}

In the following, we describe three different strategies for approximating the VEM basis functions and their gradients, and thus the spaces $\nVh{}$ and $\nabla \nVh{}$, while discussing their properties and advantages. More precisely, we consider
\begin{enumerate}
    \item the standard \textbf{NAVEM strategy}, introduced in \cite{PintoreTeora2024, PintoreTeora2025} and denoted here by the letter $\NAVEM$ to highlight the harmonic nature of NAVEM functions. Indeed, this approach approximates the VEM basis functions and their gradients as linear combinations of suitably chosen harmonic functions. As a result, Property \eqref{eq:vem:harmonicity} is exactly recovered, whereas Property \eqref{eq:vem:linearity} is only satisfied in an approximate sense, i.e. the NAVEM basis functions $\varphi^{\NAVEM}_{i} \in \con{0}{\partial E}$, but they are only approximately polynomials of degree one on each edge of the polygon.
    \item \textbf{Boundary-based NAVEM strategy (B-NAVEM)}, denoted by $\BNAVEM$. In this case, the basis functions are constructed so that Property \eqref{eq:vem:linearity} is exactly satisfied, while their harmonicity is only approximated. In particular, B-NAVEM basis functions $\varphi^{\BNAVEM}_i$ attain the exact prescribed values on the boundary of each element $E \in \Th$ and are such that $\varphi^{\BNAVEM}_i \in \con{0}{\overline{\Omega}}$.
    \item \textbf{Polynomial-based NAVEM strategy (P-NAVEM)}, denoted by $\NAPEM$. Since neither of the two previous strategies yields the inclusion $\Poly{k}{E} \subset \nVh[E]{}$ in an exact manner, we introduce a third approach, where the basis functions are again constructed to exactly satisfy Property \eqref{eq:vem:linearity}, but in this case, they are designed to directly enforce polynomial reproducibility, at least approximately.
\end{enumerate}

Each strategy relies on the same architecture for the underlying neural network, i.e.
we only consider standard fully-connected feed-forward neural networks, also known as multi-layer perceptrons \cite{goodfellow2016deep}. Given the encoding $\xx_0 \in \R^{N_0}$ of a proper dimension $N_0$ for the pair $(j,E)$, such architecture can be represented by the following formula: 
\begin{equation} 
    \begin{aligned}
        &\xx_\ell = \rho(A_\ell \xx_{\ell-1} + b_\ell), \hspace{2cm} \ell = 1,...,L-1, \\
        &{\mathcal N}(\xx_0) = A_{L} \xx_{L-1}  + b_L.
    \end{aligned}
    \label{eq:nn_formula}
\end{equation}
In this formula, the matrices $A_\ell\in\R^{N_\ell\times N_{\ell-1}}$ and the vectors $b_\ell\in\R^{N_\ell}$ contain the trainable weights of the neural network, which are optimized during the training phase, $L\in\N$ is the number of layers, $\rho:\R\rightarrow\R$ is a nonlinear scalar activation function which is applied entry-wise to the vector $A_\ell \xx_{\ell-1} + b_\ell$, and ${\mathcal N}(\xx_0) \in \R^{N_L}$ represents the neural network output.

We remark that $\mathcal N$ approximates a function mapping the vector $\xx_0\in \R^{N_0}$ to a target vector $\xx_L = \mathcal N(\xx_0)\in\R^{N_L}$. Therefore, even though $L$, $\rho$, and the intermediate values $N_{\ell}$, $\ell=1,\dots,L-1$ are model hyperparameters that can be tuned by the user, $N_0$ and $N_L$ depend on the map that the neural network approximates and cannot be modified. Thus, the different NAVEM strategies may use different values for $N_0$ and $N_L$, since the role of the underlying neural network is different in each method. In the following, we always consider the hyperbolic tangent $\rho(x) = \tanh(y)$ as activation function, and we initialize the neural network weights using the Glorot normal initialization \cite{glorot2010understanding}.

\begin{remark}\label{rem:mappingelements}
    We recall that, to improve the neural network accuracy, we always map each polygon $E$ to a reference polygon $\widehat{E}$ through the affine mapping introduced in \cite{PintoreTeora2025, Teora2024} and define all the local basis functions on $\widehat{E}$. This approach reduces geometric variability while performing input reduction, thereby improving the neural network accuracy. 
\end{remark}

\section{The NAVEM formulation}\label{sec:navem}

In this section, for the sake of completeness, we briefly summarize the standard NAVEM method. We refer the reader to \cite{PintoreTeora2025} for further details. 

Let us consider the pair $(j, E)$, we are interested in finding an approximation $\left(\varphi_{j,E}^{\rm{H}}, \qq_{j,E}^{\rm{H}}\right)$ that is cheaply computable and accurately approximates $(\varphi_{j,E}, \nabla \varphi_{j,E})$.

Let $z=x_1+ix_2$ be a complex scalar value and let $\Re(z)$ and $\Im(z)$ be its real and imaginary part. Let us consider the reference square $\refregion = [-\refedge, \refedge]^2\subset \R^2$ s.t. $E \subset \refregion$ and let $\HPoly{\refell}{\refregion}$ be the set of scaled harmonic polynomials defined over the squared region $\refregion$ (see Remark \ref{rem:mappingelements}) up to order $\refell\ge0$:
\begin{equation}
      \HPoly{\refell}{\refregion} = {\text{span}}\left\{ 1,\ \Re\left(\left(\frac{z}{\refedge}\right)^\ell\right), \Im\left(\left(\frac{z}{\refedge}\right)^\ell\right),\ \ell = 1,\dots, \refell \right\}.
      \label{eq:scaled_harmonic_polynomial}
\end{equation}
Let us further consider the auxiliary problem:
\begin{equation}\label{eq:hanging_prob}
\begin{cases}
        \Delta \tilde{\Phi} = 0 & \text{in } \Omega_{\Phi} = (-1,1)^2,\\
        \tilde{\Phi} = 1 + x_2 & \text{on } \Gamma_{\Phi, 1} = \{x_1 = 1 \text{ and } -1 \leq x_2 \leq 0\},\\
        \tilde{\Phi} = 1 - x_2 & \text{on } \Gamma_{\Phi, 2} = \{x_1 = 1 \text{ and } 0 \leq x_2 \leq 1\},\\
        \tilde{\Phi} = 0 & \text{on } \partial \Omega_{\Phi} \setminus \{\Gamma_{\Phi, 1} \cup \Gamma_{\Phi, 2}\}.
\end{cases}
\end{equation}
By solving a linear least squares problem, we look for an approximate solution of Problem \eqref{eq:hanging_prob} in the form
\begin{equation}
    \Phi(z) = \sum_{\alpha =1}^{N^{\Phi,\text{1}}} c^{\Phi,\text{1}}_{\alpha} \Re\left(\frac{d_{\alpha}}{z - z_{\alpha}}\right) + \sum_{\beta =0}^{N^{\Phi,\text{2}}} c^{\Phi,\text{2}}_{\beta} \Re\left(\left(\frac{z}{2}\right)^{\beta}\right),
    \label{eq:hanging_function}
\end{equation}
where $z_\alpha = 1 + 2 \exp\left(-4 (\sqrt{N^{\Phi,1}} - \sqrt{\alpha})\right)$, for $\alpha=1,\dots,N^{\Phi,\text{1}}$ are poles chosen to reproduce the singularity of $\Phi$ near the vertex $(1,0)$ of its domain.

We solve this problem just once on $\Omega_{\Phi}$ and then, in order to adapt this representation to a generic element $E$, we introduce three different linear transformations that map the poles of function $\Phi$ always outside $E$ and the point $(1,0)$ belonging to $\Omega_{\Phi}$ in the $(j-1)$, $j$, and $(j+1)$-th vertices of $E$, respectively. Applying these three transformations to $\Phi$, we obtain three functions, denoted by $\Phi_{j,E}^{j-1}$, $\Phi_{j,E}^{j}$ and $\Phi_{j,E}^{j+1}$, that help to improve the neural network accuracy. See Remark 3 in \cite{PintoreTeora2025} to better understand how these functions influence NAVEM accuracy.

We finally introduce the space $\approxspace[j,E]$ of harmonic functions as the space:
\begin{equation*}
     \approxspace[j,E] = \HPoly{\refell}{\refregion} \bigcup {\text{span}}\left\{\Phi_{j,E}^{j-1}, \Phi_{j,E}^{j}, \Phi_{j,E}^{j+1} \right\},
\end{equation*}
whose dimension is $\text{dim}\approxspace[j,E] = 2\refell + 4$.

To simplify the notation, we denote by $\{h_k\}_{k=1}^{2\refell+4}$, the functions spanning $\approxspace[j,E]$, namely the $2\refell+1$ basis functions defined in \eqref{eq:scaled_harmonic_polynomial} together with the three functions $\Phi_{j,E}^{j-1}$, $\Phi_{j,E}^{j}$, and $\Phi_{j,E}^{j+1}$. 
Given this approximation space $\approxspace[j,E]$, the VEM basis functions and their gradients are approximated via the standard NAVEM approach as follows: for each pair $(j, E)$,
\begin{gather}
\label{eq:navem_basis_expr} \varphi_{j,E} \approx \varphi_{j,E}^{\NAVEM} = \sum_{k=1}^{\dim \approxspace[j,E]} c_k^{\NAVEM, \varphi} h_k, \\
\nabla \varphi_{j,E} \approx \qq_{j,E}^{\NAVEM} = \sum_{k=1}^{\dim \nabla \approxspace[j,E]} c_k^{\NAVEM, \qq} \nabla h_{k+1}.
\end{gather}
The coefficients of the above linear combinations represent the outputs of two different neural networks 
\begin{equation*}
    \cc^{\NAVEM, \varphi} = \mathcal{N}^{H, \varphi}(\xx_0^{\NAVEM})\quad \text{and}\quad \cc^{\NAVEM, \qq} = \mathcal{N}^{\NAVEM, \qq}(\xx_0^{\NAVEM}),
\end{equation*}
 for a given input $\xx_0^{\NAVEM}$ that encodes the pair $(j,E)$, whose architectures are defined in \eqref{eq:nn_formula}. In particular, during the training phase, the weights of the neural network $\mathcal{N}^{\NAVEM, \qq}$ used to approximate gradients are initialized using the optimized weights of the neural network $\mathcal{N}^{\NAVEM, \varphi}$.
 \begin{remark}\label{rem:navem:gradients}
We remark that given the output of the first neural network $\{c_k^{\NAVEM, \varphi}\}_{k=1}^{\dim, \approxspace[j,E]}$ in \eqref{eq:navem_basis_expr}, then the gradients of VEM functions can be trivially approximated as
\begin{equation*}
    \nabla\varphi_{j,E} \approx \sum_{k=1}^{\dim \approxspace[j,E]} c_k^{\NAVEM, \varphi} \nabla h_{k}.
\end{equation*}
Nonetheless, we observed that the use of a second neural network to approximate the gradients improves the accuracy of the method by producing less oscillating gradients \cite{PintoreTeora2025}.
\end{remark}

 These neural networks are trained to minimize the mean squared errors  over all pairs $(j, E)$ in a suitable training dataset of the following quantities:
\begin{gather}
\label{eq:navem:loss_basis}   \epsilon_{j,E}^{\NAVEM,\varphi} = \norm[{\sob{1/2}{\partial E}}]{ \varphi_{j,E}^{\NAVEM} - \varphi_{j,E}},\\
\label{eq:navem:loss_gradient_basis}   \epsilon_{j,E}^{\NAVEM,\qq} = \norm[{\leb{2}{\partial E}}]{ \left(\qq_{j,E}^{\NAVEM}  - \nabla \varphi_{j,E}\right) \cdot {\bm{t}}},
\end{gather}
where $\tt$ denotes the unit tangent vector to the boundary $\partial E$.

As discussed in \cite{PintoreTeora2025}, these loss functions are computable since virtual basis functions are well known at the boundary of the elements. Moreover, this choice ensures good approximation properties over the entire element $E$ (see \cite{navemElasticity} for further details). We also observe that Property~\eqref{eq:vem:harmonicity} is exactly satisfied by functions defined as in \eqref{eq:navem_basis_expr}, whereas the extent to which Property~\eqref{eq:vem:linearity} is satisfied depends on the accuracy of the trained neural networks.

Note that for these NAVEM neural networks, the following requirements hold:
\begin{itemize}
    \item the output dimension must coincide with the dimension of $\approxspace[j, E]$ and of $\nabla \approxspace[j,E]$ for the two neural networks, respectively, i.e. $N_L^{\NAVEM,\varphi} = \dim \approxspace[j, E]$ and $N_L^{\NAVEM,\qq} = \dim \nabla \approxspace[j, E]$. Note that $\dim \nabla \approxspace[j, E] = \dim \approxspace[j, E] - 1$, since the kernel of $\nabla$ operator in $\approxspace[j,E]$ contains only the constant polynomial;
    \item the input dimension $N_0^{\NAVEM}$ is the same for both the neural networks and must be sufficient to properly encode the input pair $(j, E)$ into the input vector $\xx_0^{\NAVEM}$. In \cite{PintoreTeora2025}, an input reduction strategy is proposed that encodes the index of the basis function and the vertices coordinates representing the element $E$ in a vector with dimension $N_0^{\NAVEM} = 2(\Nv[E] - 1)$.
\end{itemize}

\begin{remark}\label{rem:classification}
    Since the number of vertices of an element determines the neural network architecture, in particular its input layer, the elements are implicitly grouped into classes according to their number of vertices, and a distinct neural network must be trained for each class. Nonetheless, the number of vertices per element in most meshes is typically bounded, so this assumption is not restrictive and only a limited number of neural networks is required. 
\end{remark}

\section{Enforcing continuity across adjacent elements}\label{sec:boundary_funcs}
As mentioned before, the standard NAVEM basis functions are no longer continuous across adjacent elements, since Property~\eqref{eq:vem:linearity} is enforced only approximately by minimizing the loss functions \eqref{eq:navem:loss_basis} and \eqref{eq:navem:loss_gradient_basis}. 

In this section, we show how to define an operator 
\begin{equation}
    {\mathcal B}_{j,E}:\con{0}{\overline{E}}\rightarrow \con{0}{\overline{E}}
    \label{eq:operator_B_def}
\end{equation}
that can be used to exactly enforce Property~\eqref{eq:vem:linearity} when approximating $\varphi_{j,E}$. 
The definition of the operator $\mathcal{B}_{j,E}$ is based on the definition of two auxiliary functions, namely $\psi_{j,E}$ and $\psi^0_E$. Thus, in the following, we show how to define a function $\psi_{j,E}$ that coincides with $\varphi_{j,E}$ on $\partial E$ and a bubble function $\psi^0_E$ that vanishes on $\partial E$, is strictly positive inside $E$, and has non-vanishing inward normal derivative on $\partial E$. 

Such an operator will be used in the Sections \ref{sec:bnavem} and \ref{sec:napem} to devise two alternative approaches that allow us to obtain basis functions that are continuous all across the domain $\Omega$.

In the following, the symbol $\norm{}$ denotes the Euclidean norm in $\R^2$. Moreover, we denote by $\{e_1, \dots, e_{\Nv[E]}\}$ the set of edges of $E$, where $e_i$ connects the vertices $\vv_i$ and $\vv_{i+1}$ of $E$, with $\vv_{\Nv[E] + 1} = \vv_1$, and we define $\nn_i$ as the unit outward normal vector to the edge $e_i$ with respect the polygon $E$.

\begin{figure}[!ht]
    \centering
    \begin{subfigure}{0.32\textwidth}
        \includegraphics[width=\linewidth]{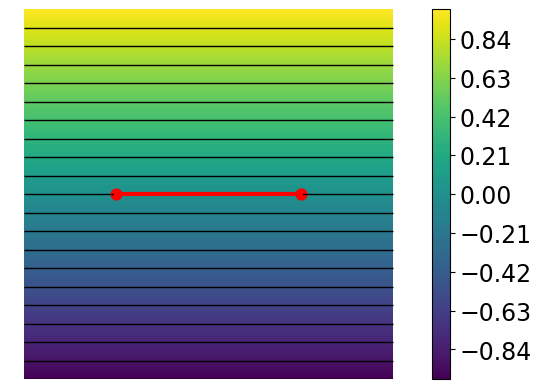}
        \caption{}
        \label{fig:bubble:distance}
    \end{subfigure}
    \begin{subfigure}{0.32\textwidth}
        \includegraphics[width=\linewidth]{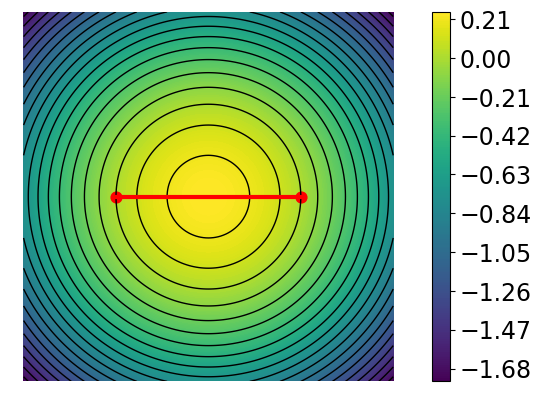}
        \caption{}
        \label{fig:bubble:trim}
    \end{subfigure}
    \begin{subfigure}{0.32\textwidth}
        \includegraphics[width=\linewidth]{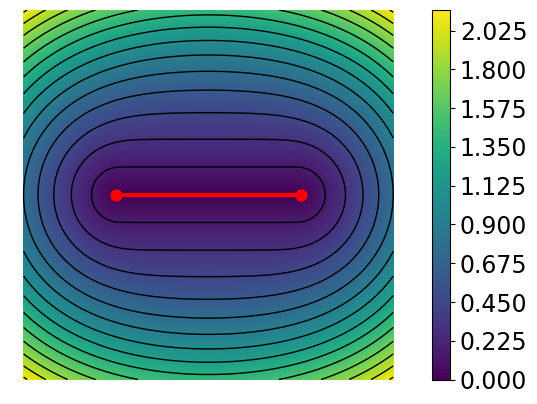}
        \caption{}
        \label{fig:bubble:adf_segment}
    \end{subfigure}
    \caption{Left: Distance function \eqref{eq:def_phi_line} to the line lying on the segment of extrema $(-0.5, 0.0)$ and $(0.5, 0.0)$. Center: trimming function defined in \eqref{eq:trim_circle} in related the the edge of extrema $(-0.5, 0.0)$ and $(0.5, 0.0)$. Right: the ADF function defined in \eqref{eq:def_phi_segment} related the the edge of extrema $(-0.5, 0.0)$ and $(0.5, 0.0)$.}
    \label{fig:bubble:edge}
\end{figure}

\subsection{A bubble function on a generic polygon}\label{sec:bubble}

Our goal is to build a bubble function $\psi_E^0$ that vanishes on $\partial E$, is strictly positive inside $E$, and has non-vanishing inward normal derivative on $\partial E$ \cite{RvachevSheiko1995RFunctionsIB}.

For this purpose, we first aim to build a function $\psi_{i,E}^0$, for each $i =1,\dots, \Nv[E]$, that is computable and strictly positive on $E$ and vanishes only on the edge $e_i$.

\begin{figure}[!ht]
    \centering
    \begin{subfigure}{0.32\textwidth}
        \includegraphics[width=\linewidth]{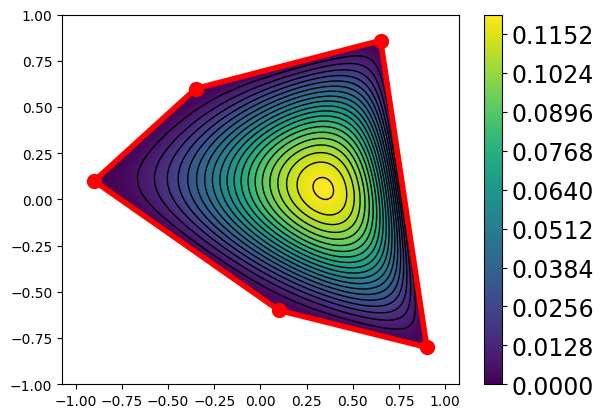}
        \caption{}
        \label{fig:bubble:convex}
    \end{subfigure}
    \begin{subfigure}{0.32\textwidth}
        \includegraphics[width=\linewidth]{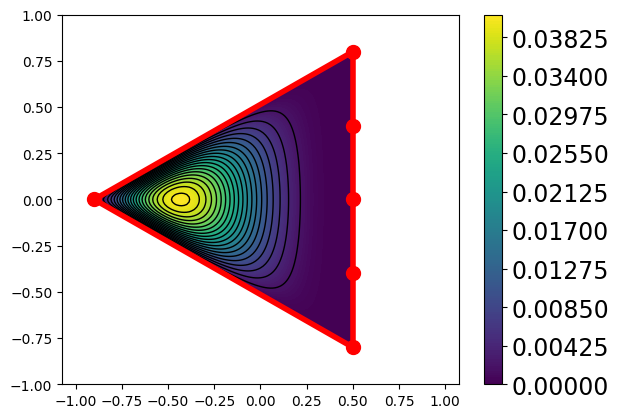}
        \caption{}
        \label{fig:bubble:hanging}
    \end{subfigure}
    \begin{subfigure}{0.32\textwidth}
        \includegraphics[width=\linewidth]{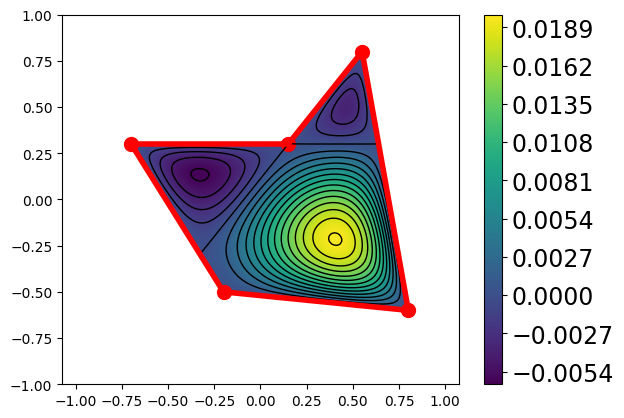}
        \caption{}
        \label{fig:bubble:concave}
    \end{subfigure}
    \begin{subfigure}{0.32\textwidth}
        \includegraphics[width=\linewidth]{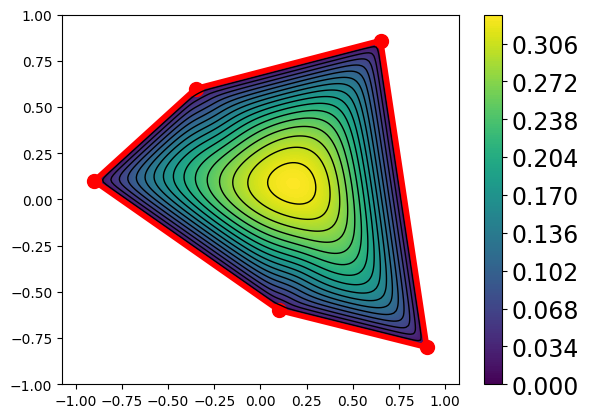}
        \caption{}
        \label{fig:adf:convex}
    \end{subfigure}
    \begin{subfigure}{0.32\textwidth}
        \includegraphics[width=\linewidth]{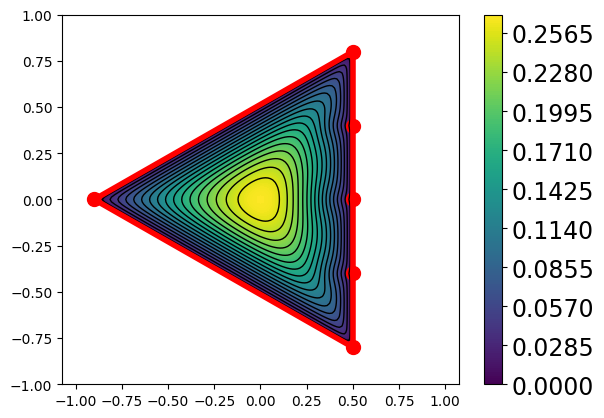}
        \caption{}
        \label{fig:adf:hanging}
    \end{subfigure}
    \begin{subfigure}{0.32\textwidth}
        \includegraphics[width=\linewidth]{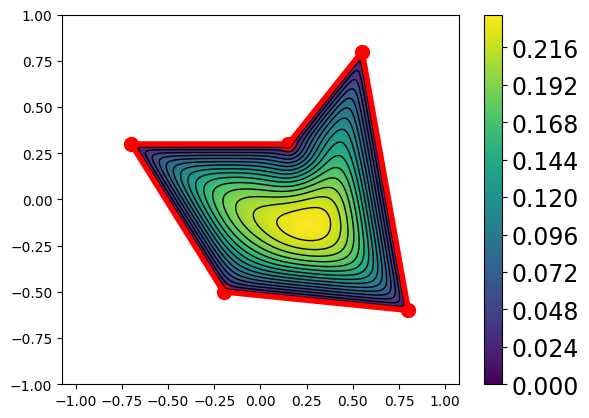}
        \caption{}
        \label{fig:adf:concave}
    \end{subfigure}
    \caption{Different alternatives to define the bubble function $\psi_{E}^0$ for a convex polygon (left), a triangle with hanging nodes (center), and a concave polygon (right). Top row refers to definition \eqref{eq:bubble_convex}, whereas the bottom row refers to definition \eqref{eq:def_phi_segment}.}
    \label{fig:bubble}
\end{figure}

To build these functions, different strategies that vary according to the shape of polygon $E$ can be pursued. If the element $E$ is convex, we can define the function $\psi_{i,E}^0$ as the \textit{signed distance function $d_i$} from the line $l_i$ where the edge $e_i$ lies, i.e. 
\begin{equation}
    \psi_{i,E}^0(\xx) := d_i(\xx) = \left(\xx - \vv_i\right)\cdot \nn_i, \quad \forall \xx \in \overline{E}.
    \label{eq:def_phi_line}
\end{equation}
which is strictly positive in $E$, and is zero on the whole $l_i$ and, in particular, on $e_i$. For an illustration see Figure \ref{fig:bubble:distance}. Given these distance functions, we can finally compute
\begin{equation}
    \psi_{E}^0(\xx) = \prod_{i=1}^{N_E^v} \psi_{i,E}^0(\xx), \quad \text{with}\quad \psi_{i,E}^0 \equiv d_i,
    \label{eq:bubble_convex}
\end{equation}
which vanishes on the boundary of $\partial E$ and is strictly positive inside $E$ (see Figure \ref{fig:bubble:convex}). This construction ensures that all the derivatives of $\psi_{E}^0$ exist and are bounded in $\overline{E}$, although $\psi_{E}^0$ may be very small in regions close to many edges of $E$ and has zero normal derivative. For an example, see Figure \ref{fig:bubble:hanging}.

Moreover, if the element $E$ is concave, the line $l_i$ where the signed distance function $d_i$ is zero could cross the polygon, making $\psi_{i,E}^0$ no longer be strictly positive inside $E$ as requested (see Figure \ref{fig:bubble:concave}). In these cases, an example of a function that vanishes only on $e_i$ is represented by the \textit{Approximate Distance Function} (ADF) defined in \cite{Pintore2023}. More specifically, let us define the trimming function $t_i^{\rm{c}}$ as
\begin{equation}
    t_i^{\rm{c}}(\xx) = \frac {1}{\vert e_i\vert} \left[\left(\frac {\vert e_i \vert}{2}\right)^2 - \norm{\xx - \frac{\vv_i+\vv_{i+1}}{2}}^2\right], \quad \forall \xx \in \overline{E}.
    \label{eq:trim_circle}
\end{equation}
Note that $t_i^{\rm{c}}\ge 0$ defines a circle of center $\dfrac{\vv_i+\vv_{i+1}}{2}$ as shown in Figure \ref{fig:bubble:trim}. Now, we can define $\psi_{i,E}^0$ as the ADF function related to $e_i$, i.e.
\begin{equation}
    \psi_{i,E}^0 = \sqrt{d_i^2(\xx) + \left(\frac{\sqrt{(t_i^{\rm{c}}(\xx))^2+d_i^4(\xx)}-t_i^{\rm{c}}(\xx)}{2}\right)^2}, \quad \forall \xx \in \overline E.
    \label{eq:def_phi_segment}
\end{equation}
By defining $\psi_{i,E}^0$ as the ADF related to $e_i$, we obtain a function that is zero only on the edge $e_i$ and is strictly positive elsewhere. For an illustration of such a function, see Figure \ref{fig:bubble:adf_segment}.

Finally, the bubble function $\psi_E^0$ can be chosen as the ADF to $\partial E$, normalized up to order $m = 2$, i.e. as
\begin{equation}
    \psi_E^0(\xx) = \dfrac{1}{\sqrt[m]{{\mathlarger{\mathlarger{\sum}}}_{i=1}^{\Nv[E]} \dfrac{1}{\left(\psi_{i,E}^0(\xx)\right)^{m}}} }. 
    \label{eq:adf_boundary}
\end{equation}
With this definition, $\psi_E^0$ vanishes at the boundary, is strictly positive inside $E$, and has a unitary inward normal derivative at the boundary, independently of the shape of the polygon $E$. See Figures \ref{fig:adf:convex}, \ref{fig:adf:hanging}, and \ref{fig:adf:concave} for an illustration of the ADF function for different kinds of polygons.
Moreover, it is $\con{2}{}$ in all the points away from the boundary $\partial E$, while its Laplacian blows up at the vertices $\vv_i$ (i.e. at the points belonging to $l_i$ such that $t_i^{\rm{c}}(\xx)=0$). 

Since the combination \eqref{eq:def_phi_segment}-\eqref{eq:adf_boundary} is applicable to all classes of polygons and has been numerically observed to provide better accuracy, it is therefore employed regardless of the specific polygon under consideration, i.e. for all the convex and concave polygons, as well as for polygons characterized by hanging nodes.

\subsection{The transfinite interpolation of $\varphi_{j,E}$}\label{sec:boundary_funz}

Let us consider a point $\xx \in \overline{E}$ and let $\Pi_{e_i} \xx$ be the orthogonal projection of $\xx$ onto the line $l_i$ where the edge $e_i$ lies. For any $\zz \in l_i$, let $s_i(\zz)$ be the curvilinear coordinate of $\zz$ with respect to the origin $\vv_i$ and the direction $(\vv_{i+1}-\vv_i)$, i.e. $s_i$ is a linear function such that $s_i(\vv_i)=0$ and $s_i(\vv_{i+1})=1$, whose expression is:
\begin{equation*}
    s_i(\zz) = \frac{(\zz - \vv_i) \cdot (\vv_{i+1} - \vv_i)}{\norm{\vv_{i+1} - \vv_i}^2}\quad \forall \zz \in l_i.
\end{equation*}
Let $\xx_{i,\perp} := \xx - \Pi_{e_i} \xx$ be the component of $\xx$ orthogonal to the vector $\vv_{i+1} - \vv_i$. The scalar product $(\xx - \vv_i) \cdot (\vv_{i+1} - \vv_i)$ can be written as:
\begin{equation*}
    \begin{aligned}
        (\xx - \vv_i) \cdot (\vv_{i+1} - \vv_i) &= (\xx_{i,\perp} + \Pi_{e_i} \xx - \vv_i) \cdot (\vv_{i+1} - \vv_i) \\
        &= (\Pi_{e_i} \xx - \vv_i) \cdot (\vv_{i+1} - \vv_i) + \xx_{i,\perp} \cdot (\vv_{i+1} - \vv_i) \\
        &= (\Pi_{e_i} \xx - \vv_i) \cdot (\vv_{i+1} - \vv_i).
    \end{aligned}
\end{equation*}
As a consequence, we can define the function $s_i : \overline{E} \to \R$ as
\begin{equation*}
    s_i(\Pi_{e_i}\xx) = \frac{(\xx - \vv_i) \cdot (\vv_{i+1} - \vv_i)}{\norm{\vv_{i+1} - \vv_i}^2} \quad \forall\xx\in\overline{E}.
\end{equation*}

For any $i,j=1,\dots,\Nv[E]$ and for any $\xx\in \overline E$, we introduce the functions
\begin{equation*}
    \psi_{i,j,E}(\xx) = \begin{cases}
        s_i(\Pi_{e_i} \xx) & \text{if } j = i-1,\\
        1 - s_i(\Pi_{e_i} \xx) & \text{if } j = i,\\
        0 & \text{otherwise}.
    \end{cases}
\end{equation*}
These functions are defined everywhere in $\overline E$ and are known in a closed form. Moreover, they are linear polynomials on each edge and satisfy 
\begin{equation*}
    \psi_{i,j,E}(\vv_j) = 1 \quad\text{and}\quad \psi_{i,j,E}(\vv_k) = 0 \quad \forall k \neq j.
\end{equation*}
Thus, it can be easily proved that $\psi_{i,j,E}=\varphi_{j,E}$ on the edge $e_i$. Then, the transfinite interpolation $\psi_{j,E}$ \cite{Sukumar2021} of the function $\varphi_{j,E}$ is defined as
\begin{equation}
    \psi_{j,E}(\xx) = \sum_{j=1}^{\Nv[E]} \omega_{j,E}(\xx) \psi_{i,j,E}(\xx),
    \label{eq:transfite_interpolant}
\end{equation}
where
\begin{equation*}
    w_{j,E}(\xx) = \dfrac{\prod_{i=1; i \ne j}^{\Nv[E]} \psi_{i,E}^0(\xx)}{\sum_{r=1}^{\Nv[E]} \prod_{i=1;i\ne r}^{\Nv[E]} \psi_{i,E}^0(\xx)}.
\end{equation*}
and the functions $\psi_{i,E}^0$ are defined in \eqref{eq:def_phi_segment}.
Given the properties of these functions, it is immediate to check that
\begin{equation*}
    \psi_{j,E} = \varphi_{j, E} \quad \text{on } \partial E.
\end{equation*}

\subsection{The operator $\mathcal{B}_{j,E}$}

Given the bubble function $\psi_E^0$, defined in \eqref{eq:def_phi_segment}-\eqref{eq:adf_boundary} and vanishing on $\partial E$, and the transfinite interpolation $\psi_{j,E}$ of $\varphi_{j,E}$, introduced in \eqref{eq:transfite_interpolant} and satisfying $\psi_{j,E} = \varphi_{j,E}$ on $\partial E$, we define the operator ${\mathcal B}_{j,E}$, introduced in \eqref{eq:operator_B_def}, as
\begin{equation}
    \left({\mathcal B}_{j,E}(v)\right)(\xx) := \psi_E^0(\xx) v(\xx) + \psi_{j,E}(\xx),\quad \forall \xx \in E,\quad \forall v \in \con{0}{\overline{E}}.
    \label{eq:bndry_operator}
\end{equation}
By construction, it satisfies
\begin{equation*}
    \left({\mathcal B}_{j,E}(v)\right)(\xx) = \varphi_{j,E}(\xx), \quad \forall \xx\in \partial E, \quad \forall v \in \con{0}{\overline{E}}.
\end{equation*}

If the neural approximation $\varphi_{j,E}^{\NN}$ is defined as the application of ${\mathcal B}_{j,E}$ to the neural network output $\mathcal{N}(\xx_0)$, i.e.,
\begin{equation*}
    \varphi^{\NN}_{j,E} = {\mathcal{B}_{j,E}}(\mathcal{N}(\xx_0)),
\end{equation*}
then $\varphi_{j,E}^{\NN}$ exactly matches $\varphi_{j,E}$ on the boundary $\partial E$. Since the bubble function $\psi_E^0$ vanishes on $\partial E$, the neural network output $\mathcal{N}(\xx_0)$ influences the approximation only in the interior of the element $E$. Consequently, $\mathcal{N}(\xx_0)$ can be trained to control the interior behaviour of $\varphi_{j,E}^{\NN}$ to minimize a prescribed cost functional.

\begin{remark}
In the case of convex polygonal elements, one could employ the standard polynomial bubble function defined in \eqref{eq:bubble_convex}. However, since this bubble may attain very small values in large portions of the element, particularly in points $z$ that are close to many edges, the contribution of the neural network output 
$\mathcal{N}(\xx_0)$ is strongly damped. As a consequence, learning an effective interior correction becomes difficult, because the product $\psi_E^0(\zz) \mathcal{N}(\xx_0)$ has only a marginal influence on the final neural approximation $\left({\mathcal B}_{j,E}(\mathcal{N}(\xx_0))\right)(\zz)$.
\end{remark}

\section{The B-NAVEM formulation}\label{sec:bnavem}

The virtual basis functions $\{\varphi_{j,E}\}_{j=1}^{\Nv[E]}$ are defined as the solutions of the following local Laplace problems 
\begin{equation}
    \begin{cases}
        \Delta \varphi_{j,E} = 0 & \text{in } E,\\
        \varphi_{j,E} = \psi_{j,E} & \text{on } \partial E.
    \end{cases}
    \label{eq:local_PDE_problem}
\end{equation}

As discussed in the previous section, the map associating each input pair $(j,E)$ with the approximation function $\varphi_{j,E}^{\NN}$ is highly nonlinear. For this reason, Problem~\eqref{eq:local_PDE_problem} can be efficiently addressed using neural networks of the form~\eqref{eq:nn_formula}. Among the most widely used neural networks-based PDE solvers, Physics-Informed Neural Networks (PINN), originally introduced in~\cite{PINN}, have recently obtained a lot of attention. The core idea of the B-NAVEM method is to employ a PINN to learn the map~\eqref{eq:nn_map}.

Unlike the NAVEM approach, where the neural network output represents the vector of coefficients with respect to the harmonic basis for $\approxspace[j,E]$, in PINNs the output directly represents the value of the function at a given point $\zz$ inside the domain $E$. Consequently, the PINN effectively learns the following nonlinear map:
\begin{equation}
    (\zz, j, E) \mapsto \varphi^{\NN}_{j,E}(\zz).
    \label{eq:nn_map_pinn}
\end{equation}
Thus, the input and output dimensions in the B-NAVEM method are $N_0^{\BNAVEM} = 2 \Nv[E]$ and $N_L^{\BNAVEM} = 1$, respectively. In particular, the input vector $\xx_0^{\BNAVEM}$ consists of the concatenation of the evaluation point $\zz \in \R^2$ and the NAVEM encoding $\xx_0^{\NAVEM} \in \R^{2(\Nv[E] - 1)}$ of the pair $(j,E)$. More precisely, we set
\begin{equation*}
    \xx_0^{\BNAVEM} = [\zz;\, \xx_0^{\NAVEM}] \in \R^{N_0^{\BNAVEM} }.
\end{equation*}

Given the operator $\mathcal{B}_{j,E}$ defined in~\eqref{eq:bndry_operator}, Dirichlet boundary conditions are enforced in the PINN framework by adopting the technique described in~\cite{Sukumar2021, Pintore2023}. Specifically, this is achieved by suitably modifying the neural network output. For an input $\xx_0^{\BNAVEM}$ encoding the triplet $(\zz, j, E)$, the B-NAVEM basis functions are defined as
\begin{equation*}
    \varphi_{j,E}^{\BNAVEM}(\zz)
    = \mathcal{B}_{j,E}\bigl(\mathcal{N}^{\BNAVEM}(\xx_0^{\BNAVEM})\bigr)(\zz)
    = \psi_E^0(\zz)\,\mathcal{N}^{\BNAVEM}(\xx_0^{\BNAVEM}) + \psi_{j,E}(\zz),
\end{equation*}
where $\mathcal{N}^{\BNAVEM}(\xx_0^{\BNAVEM}) \in \R$ denotes the scalar output of the PINN.

Since, for any choice of the trainable weights, the function $\varphi_{j,E}^{\BNAVEM}$ automatically satisfies Property~\eqref{eq:vem:linearity}, the training procedure only needs to enforce Property~\eqref{eq:vem:harmonicity}. To this end, we define the loss term as the PDE residual
\begin{equation}
    \epsilon_{j,E}^{\BNAVEM}
    = \norm[\leb{2}{\Omega}]{\Delta \varphi_{j,E}^{\BNAVEM}}.
    \label{eq:bnavem:loss}
\end{equation}

In the B-NAVEM approach, the Laplacian $\Delta \varphi_{j,E}^{\BNAVEM}(\zz)$, as well as the gradient $\nabla \varphi_{j,E}^{\BNAVEM}(\zz)$, can be efficiently computed via automatic differentiation~\cite{baydin2018automatic}. Accordingly, the discrete B-NAVEM gradient is defined as
\begin{equation*}
    \qq^{\BNAVEM}_{j,E}(\zz) = \nabla \varphi_{j,E}^{\BNAVEM}(\zz).
\end{equation*}

\section{The P-NAVEM formulation}\label{sec:napem}

As shown in the previous sections, the NAVEM approach allows us to exactly enforce the Property~\eqref{eq:vem:harmonicity}, while Property~\eqref{eq:vem:linearity} is satisfied only up to the NAVEM accuracy. Conversely, in the B-NAVEM approach, Property~\eqref{eq:vem:linearity} is exactly enforced, whereas Property~\eqref{eq:vem:harmonicity} depends on the approximation properties of the neural network. In both cases, however, we can conclude that
\begin{equation}
    \Poly{k}{E} \not\subset \nVh[E]{\ast} \quad \forall \ast \in \{\NAVEM, \BNAVEM\}, 
\end{equation}
where $ \nVh[E]{\NAVEM}$ and $ \nVh[E]{\BNAVEM}$ represent the local H-NAVEM and B-NAVEM spaces, respectively, and are defined as dictated in \eqref{eq:local_navem_space}.
More precisely, the capability of NAVEM or B-NAVEM to reproduce polynomial functions depends on the accuracy of the underlying neural networks. The lack of exact polynomial inclusion may therefore limit the convergence properties of both methods.

The theoretical VEM analysis suggests introducing a linearly independent set of basis functions $\{\varphi_{j,E}^{\NN}\}_{j=1}^{\Nv[E]}$ satisfying Property~\eqref{eq:vem:linearity}, which guarantees $\con{0}{}$-conformity, together with the two following additional properties: 
\begin{enumerate}[label=\textbf{P.\arabic*}]
    \item \label{prop:partitionofunity}\textit{Partition of unity property}:
    \begin{equation}\label{eq:partition_unity}
        \sum_{j=1}^{\Nv[E]} \varphi_{j,E}^{\NN} (\xx) = 1 \quad \forall \xx \in \overline{E},
    \end{equation}
    \item \label{prop:linearreproduction}\textit{Linear reproduction property}:
    \begin{equation}\label{eq:x_approx}
        \sum_{j=1}^{\Nv[E]}\, (\vv_j)_1\,\, \varphi_{j,E}^{\NN}(\xx) = x_1, \quad\quad \forall \xx=(x_1,x_2) \in \overline{E},
    \end{equation}
    \begin{equation}\label{eq:y_approx}
        \sum_{j=1}^{\Nv[E]}\,(\vv_j)_2\,\, \varphi_{j,E}^{\NN}(\xx) = x_2, \quad\quad \forall \xx=(x_1,x_2) \in \overline{E},
    \end{equation}
    where $((\vv_j)_1, (\vv_j)_2)$ are the coordinates of the $j$-th vertex of $E$.
\end{enumerate}
The two Properties \ref{prop:partitionofunity} and \ref{prop:linearreproduction} imply the exact linear polynomial reproducibility. This last condition is sufficient to ensure the desired convergence rates of the method, without requiring the basis functions to be harmonic as in the virtual element framewrok. 

Motivated by these considerations, we aim to construct a discrete space whose basis functions satisfy, as accurately as possible, Properties~\eqref{eq:vem:linearity}, \ref{prop:partitionofunity}, and \ref{prop:linearreproduction}. We refer to this approach as P-NAVEM, emphasizing its focus on polynomial reproducibility.

Property~\eqref{eq:vem:linearity} can be imposed exactly by designing a neural network that learns the nonlinear mapping~\eqref{eq:nn_map_pinn} and by applying the operator $\mathcal{B}_{j,E}$ to the P-NAVEM network output $\mathcal{N}^{\NAPEM}(\xx_0^{\NAPEM})$, as done in the B-NAVEM approach. Since the P-NAVEM network input encodes the same geometric information as in B-NAVEM, we have $\xx_0^{\NAPEM}=\xx_0^{\BNAVEM}$ and $N_0^{\NAPEM}=N_0^{\BNAVEM}$. Moreover, as in B-NAVEM, we set $N_L^{\NAPEM} = 1$ and define
\begin{equation}
    \varphi_{E}^{\NAPEM}(\zz)
    = \mathcal{B}_{j,E}\bigl(\mathcal{N}^{\NAPEM}(\xx_0^{\NAPEM})\bigr)(\zz)
    = \psi_E^0(\zz)\,\mathcal{N}^{\NAPEM}(\xx_0^{\NAPEM}) + \psi_{j,E}(\zz).
    \label{eq:napem:basis_functions}
\end{equation}

The key difference between B-NAVEM and P-NAVEM concerns the loss function considered during the training phase. In the P-NAVEM approach, we aim to enforce the two properties \ref{prop:partitionofunity} and \ref{prop:linearreproduction} and thus, to minimize the following quantities
\begin{equation*}
    \epsilon_{j,E}^{p_i} = \norm[\leb{2}{E}]{\sum_{j=1}^{\Nv[E]} p_i(\vv_j)\varphi^{\NAPEM}_{j,E} - p_i }, \quad \forall i\in\{0,1,2\},
\end{equation*}
where, for all $\xx \in \overline{E}$,
\begin{equation*}
p_0(\xx) = 1, \quad p_1(\xx) = x_1, \quad p_2(\xx) = x_2.
\end{equation*}
We observe that, independently of the method $\ast \in \{\NAVEM, \BNAVEM, \NAPEM\}$, Property~\ref{prop:partitionofunity} can always be enforced exactly by learning (for instance) only the first $\Nv[E] - 1$ basis functions and then defining the last one as
\begin{equation*}
    \varphi_{\Nv[E],E}^{\ast} = 1 - \sum_{j=1}^{\Nv[E] - 1} \varphi_{\Nv[E],E}^{\ast}.
\end{equation*}

To approximately enforce \ref{prop:linearreproduction}, we can thus train neural network ${\mathcal N}^{\NAPEM}$ to minimize just
\begin{equation}
    \epsilon_{E}^{\NAPEM, \varphi} = \epsilon_{E}^{p_1} + \epsilon_{E}^{p_2}.
    \label{eq:napem:loss_basis}
\end{equation}
However, as already observed for the NAVEM method, training the neural network by directly learning the basis functions $\varphi_{j,E}^{\NAPEM}$ through a loss contribution that depends exclusively on the basis functions and does not involve their gradients, namely~\eqref{eq:napem:loss_basis}, and subsequently approximating the gradients as $\nabla \varphi_{j,E} \approx \nabla \varphi_{j,E}^{\NAPEM}$, typically produces smooth basis functions but highly oscillatory gradients (see Remark~\ref{rem:navem:gradients}).

To mitigate this issue, rather than minimizing a loss functional that depends solely on the basis functions, we introduce a loss contribution that directly involves their gradients, thereby controlling their oscillatory behaviour:
\begin{equation}
\epsilon_{E}^{\NAPEM, \qq} = \epsilon_{E}^{\nabla p_1} + \epsilon_{E}^{\nabla p_2},
    \label{eq:napem:loss_gradients}
\end{equation}
where
\begin{equation*}
    \epsilon_{E}^{\nabla p_i} = \norm[\leb{2}{E}]{\sum_{j=1}^{N_E^v} p_i(\vv_j)\nabla \varphi^{\NAPEM}_{j,E} - \nabla p_i }, \quad \forall i\in\{0,1,2\}.
\end{equation*}

Moreover, in the NAVEM method we employ two neural networks to specifically approximate the basis functions and their gradients, respectively, since it has been observed to be the best alternative in terms of accuracy \cite{PintoreTeora2025}. In the P-NAVEM method, instead, we can choose to employ only one neural network that minimizes the loss term \eqref{eq:napem:loss_gradients} and that approximates in a given point $\zz \in E$, the basis functions as $\varphi_{j,E}^{\NAPEM}$ as in \eqref{eq:napem:basis_functions} and their gradients as
\begin{equation}
    \nabla \varphi_{j,E}(\zz) \approx \qq_{j,E}^{\NAPEM}(\zz) = \nabla \varphi_{j,E}^{\NAPEM}(\zz).
    \label{eq:napem:grad_basis_functions} 
\end{equation}
This can be done since the operator ${\mathcal B}_{j,E}$ forces the basis functions $\varphi_{j,E}^{\NAPEM}= {\mathcal B}_{j,E}({\mathcal N}^{\NAPEM})$ to be piece-wise linear on the boundary of the element, Lagrangian in the vertices, and, as a consequence, linearly independent.

In contrast, within the NAVEM framework, it is necessary to employ two separate neural networks for learning virtual functions and their gradients. Indeed, if only the network minimizing the gradient-based loss~\eqref{eq:navem:loss_gradient_basis} is used, the resulting gradients are less oscillatory, but the reconstructed basis functions lose the correct VEM scaling, since the constant component of their linear combination is not directly controlled. 

 

\section{Numerical results}\label{sec:numeri_experiements}

In this section, we present a series of numerical experiments aimed at comparing the proposed techniques for approximating basis functions on general polygonal meshes. 

Specifically, three numerical tests are conducted. The first experiment evaluates the training and generalization performance of neural networks across different classes of polygonal elements. The second experiment compares the accuracy and computational efficiency of the corresponding numerical methods when applied to solve a diffusion–advection–reaction problem. Finally, the third experiment highlights the advantages of the proposed approaches over the standard virtual element method, particularly in terms of eliminating projection and stabilization operators. For details on the virtual element discretization employed in this manuscript, we refer the reader to \cite{LBe16, LBe14}.

\subsection{Test 1: The neural networks and the training phase}
As highlighted in Remark \ref{rem:classification}, we have to train several neural networks to predict the basis functions, one for each class of polygons we consider. Since neural networks require a fixed input dimension and this dimension depends on the number of vertices, we have to train neural networks for each distinct class corresponding to a value of $\Nv[E]\geq 4$ appearing in our test mesh. The case $\Nv[E]=3$ does not require a neural network, as the virtual element method for the lowest order reduces to the finite element method, whose basis functions are polynomials known in closed form. Therefore, the piecewise linear finite element basis functions are also used in the NAVEM framework in the case of triangles.

\begin{figure}[!ht]
\centering
\begin{subfigure}{0.33\textwidth}
    \includegraphics[width=\textwidth]{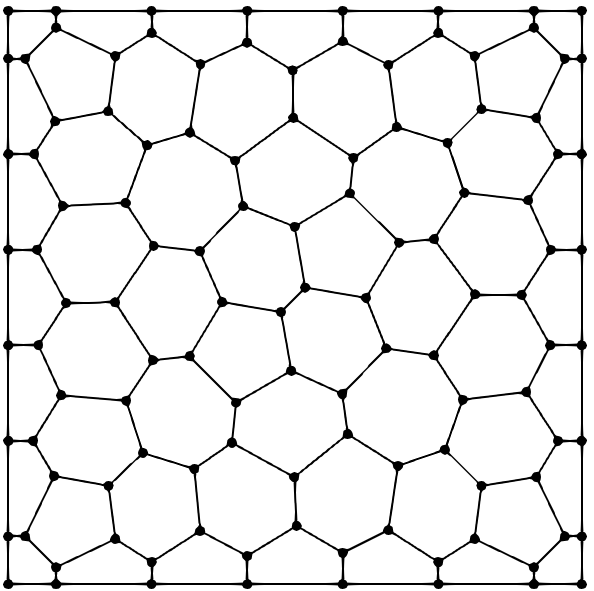}
    \caption{}
    \label{fig:voronoi_mesh}
\end{subfigure}
\begin{subfigure}{0.33\textwidth}
    \includegraphics[width=\textwidth]{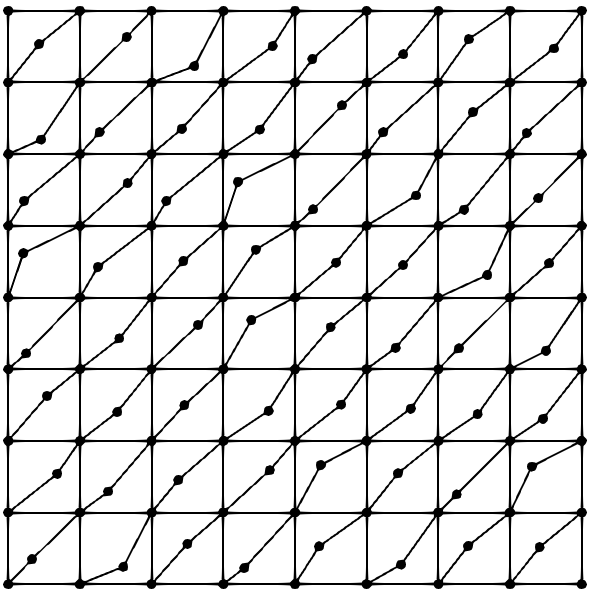}
    \caption{}
    \label{fig:rand_dist_concave_mesh}
\end{subfigure}
\caption{Tests 1, 2, and 3: Second mesh in the two considered families of meshes used in the test phase. Left: Voronoi. Right: Convex-Concave.}
\label{fig:meshes}
\end{figure}

In all the numerical tests, we consider two sets of four meshes: a family of Voronoi meshes and a family of quadrilateral meshes that comprises both concave and convex quadrilaterals. The former is generated through the MATLAB library mVEM \cite{mVEM}, while the latter is generated by randomly perturbing the vertices coordinates of a related family of structured convex-concave meshes. The second mesh for each family is shown in Figure \ref{fig:meshes}. As it can be observed from Figure \ref{fig:rand_dist_concave_mesh}, the strategy followed to build the Convex-Concave family $\{\mathcal{T}_{h,i}^Q\}_{i=1}^4$ can produce hanging nodes. In these meshes, the number of vertices is at most $7$. To further improve the accuracy of the neural network approximations, we train separate networks for convex and concave quadrilaterals. In particular, we adopt the model trained for concave elements in the case of a quadrilateral that has the shape of a triangle and has one hanging node. Consequently, for each model, we train five neural networks corresponding to the following classes of polygons:
\begin{enumerate}
    \item convex quadrilaterals;
    \item concave quadrilaterals;
    \item convex pentagons;
    \item convex hexagons;
    \item convex heptagons.
\end{enumerate}

For the first two classes, thanks to the limited geometric variability of quadrilateral elements,  we generate synthetic datasets consisting of $1000$ randomly generated convex and concave quadrilaterals, respectively, by using the Python library \href{https://pypi.org/project/polygenerator/}{polygenerator}. For the remaining classes, instead, we sample elements with the prescribed number of vertices from training Voronoi meshes generated using the same algorithm adopted for the test meshes. Specifically, we obtain $1000$ elements for convex pentagons and hexagons, and about $300$ elements for heptagons, since they are rarer in Voronoi meshes. This strategy ensures that the training polygons closely resemble those encountered in the test phase, thereby enabling the networks to generalize effectively to unseen elements. We remark that these datasets (and the corresponding NAVEM networks) are the same employed to perform the simulation in \cite{navemElasticity} on a different test family of Voronoi meshes.

\begin{algorithm}
\caption{Algorithm designed to select points inside the triangle of vertices $A=(0,0)$, $B=(1,0)$, $C=(0,1)$ to approximate loss integrals.}
\label{alg:points_on_triangle}
\KwData{$N \ge 0$,  $\varepsilon>0$}
\KwResult{List of $(N+1)(N+2)/2$ points}
$\rm{points} \gets [\,]$\;
\smallskip
\For{$0 \le i \le N$}{
    \smallskip
    \For{$0 \le j \le N-i$}{
        \smallskip
        $k \gets N-i-j$\;
        \smallskip
        $x_0 \gets \dfrac{i+0.5}{N+1.5}$\;
        \smallskip
        $y_0 \gets \dfrac{j+0.5}{N+1.5}$\;
        \smallskip
        $z_0 \gets \dfrac{k+0.5}{N+1.5}$\;
        \smallskip

        $z \gets 1 - \cos\left(\dfrac{\pi}{2}   z_0\right)$ \Comment*[r]{biased version of $z$ to accumulate points close to edge $\overline{BC}$}
        \smallskip
        
        $s \gets \dfrac{1 - z}{1 - z_0 + \varepsilon}$ \Comment*[r]{scaling factor}
        \smallskip

        $\rm{new\_point} \gets \left(s\,x_0,\, s\,y_0\right)$ \Comment*[r]{define the new point}
        \smallskip
        $\rm{points}.\rm{insert}(\rm{new\_point})$ \Comment*[r]{add the new point to the list}
      }
    }
\Return{$\rm{points}$}
\end{algorithm}

\begin{figure}[!ht]
\centering
\begin{subfigure}{0.33\textwidth}
    \includegraphics[angle=-90, width=\textwidth]{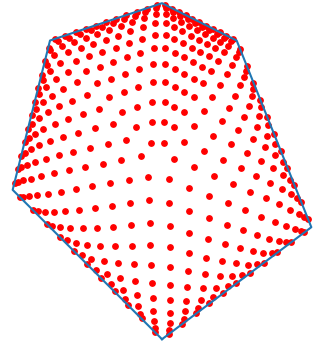}
    \caption{}
    \label{fig:quad_pts_convex}
\end{subfigure}\hspace{70pt}
\begin{subfigure}{0.33\textwidth}
    \includegraphics[angle=90, width=\textwidth]{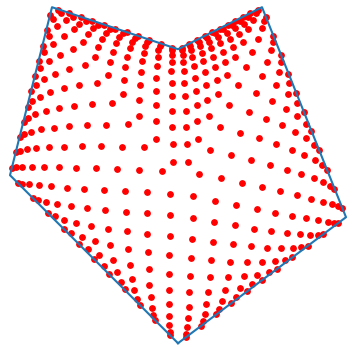}
    \caption{}
    \label{fig:quad_pts_concave}
\end{subfigure}
\caption{Points distribution generated by using the Algorithm \ref{alg:points_on_triangle} on general polygons.}
\label{fig:quad_pts}
\end{figure}

For each class of polygons, we train our three models, namely NAVEM, B-NAVEM, and P-NAVEM, by minimizing the related losses on the training data. The training procedure requires the numerical evaluation of the integrals appearing in the loss definitions. 

In the case of the NAVEM method, the integrals defining the loss functions \eqref{eq:navem:loss_basis} and \eqref{eq:navem:loss_gradient_basis} are approximated using a set of $50$ quadrature points on each edge of the polygon, distributed exponentially toward the polygon vertices.

For the B-NAVEM and P-NAVEM methods, the loss integrals \eqref{eq:bnavem:loss} and \eqref{eq:napem:loss_gradients} require an appropriate selection of points inside each element. To this end, we first construct a simple triangulation of the polygon by connecting all its vertices to a single interior point. The existence of such a point is guaranteed under the standard VEM assumption that elements are star-shaped \cite{LBe13, LBe17}. On each resulting sub-triangle, we then apply Algorithm~\ref{alg:points_on_triangle} with $N=10$ to generate interior sampling points. Although the algorithm is formulated for the reference triangle with vertices $A=(0,0)$, $B=(1,0)$, and $C=(0,1)$, points in any physical triangle are readily obtained via the standard finite element mapping from the reference configuration. The resulting point distribution within a generic polygon tends to cluster near the element edges, as illustrated in Figure~\ref{fig:quad_pts_convex}. Moreover, the point density naturally increases close to the boundary and in sub-triangles associated with short edges or near re-entrant corners (see Figure~\ref{fig:quad_pts_concave}). This feature is particularly advantageous, since the target functions typically exhibit larger gradients in these regions; a denser sampling therefore improves the accuracy of the neural network approximation in these areas.

\begin{figure}[!ht]
\centering
\begin{subfigure}{0.48\textwidth}
    \includegraphics[width=\textwidth]{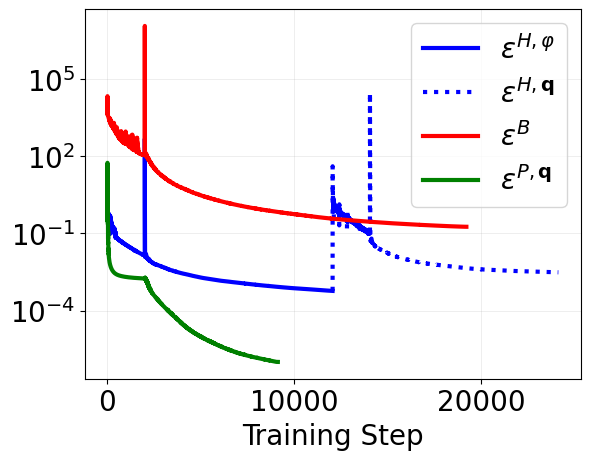}
    \caption{}
    \label{fig:loss_vs_epochs}
\end{subfigure}
\begin{subfigure}{0.48\textwidth}
    \includegraphics[width=\textwidth]{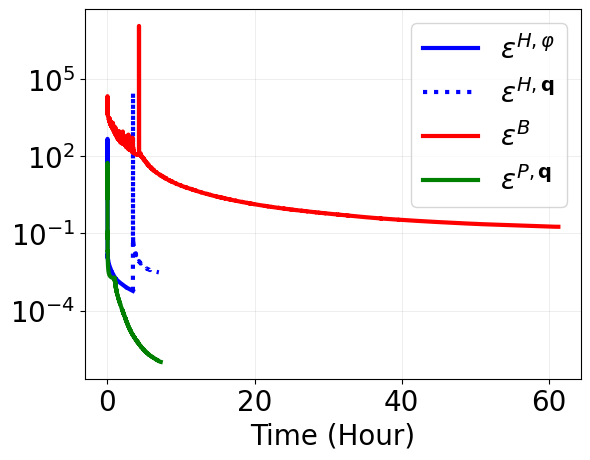}
    \caption{}
    \label{fig:loss_vs_time}
\end{subfigure}
\caption{Test 1: Behaviour of loss values as the number of training steps increases (left) and time increases (right) during the training phase.}
\label{fig:losses_training}
\end{figure}

To train each network, we use 2000 epochs of the Adam optimizer \cite{kingma2014adam} with an exponentially decaying learning rate from $10^{-2}$ to $10^{-3}$, and up to 10000 epochs of the BFGS optimizer \cite{Battiti1990BFGS}. The evolution of the loss functions, either with respect to the number of weight updates and the elapsed training time, is reported in Figure \ref{fig:losses_training} for the case of convex quadrilaterals. We do not report the other training plots since the behaviours related to the other polygon classes are very similar to the case of the convex quadrilateral. From these results, we can observe that the NAVEM method requires a larger number of epochs to converge, but that each step is very efficient in terms of time. In contrast, the B-NAVEM method is significantly more expensive, with a total training time approximately eight times larger than that of the other two approaches.

To better explain these differences regarding the computational times, we first observe that all matrices and vectors that do not depend on the neural network parameter are precomputed in a pre-processing step. For example, the harmonic functions in $\approxspace[j,E]$ and their derivatives for NAVEM, as well as the ADF functions \eqref{eq:adf_boundary}, the transfinite interpolant \eqref{eq:transfite_interpolant}, and their derivatives for B-NAVEM and P-NAVEM, are evaluated only once at quadrature points in this pre-processing stage and subsequently reused multiple times throughout the training process to compute the losses. Since the pre-processing time is negligible with respect to training time, the pre-processing time is not accounted for in the reported time.  

Moreover, during the training process, only the neural networks and their derivatives are evaluated. Consequently, the differences in the reported computational time and effort among the three methods depend mainly on the order of neural network derivatives required. In particular, the NAVEM method requires only the evaluation of the network output, the P-NAVEM method requires both the output and its gradient, whereas the B-NAVEM method additionally requires the computation of the network Laplacian at each epoch, which is considerably more expensive. Thus, the difference in the computational times reported in Figure \ref{fig:loss_vs_time} strongly depends on the computation of the network Laplacian.

We anticipate that this disparity affects only the training phase, as will be shown in the next numerical experiments. Indeed, during the test phase, i.e. when the trained networks are employed to solve the discrete Problem~\eqref{eq:discrete_var_prob}, the computation of the Laplacian is no longer required for the B-NAVEM method, rendering the differences in computational time between B-NAVEM and P-NAVEM negligible.

\begin{remark}
The B-NAVEM network ${\mathcal N}^{\BNAVEM}$ is trained as a PINN on a parameterized polygonal domain to solve the local Laplace problem~\eqref{eq:local_PDE_problem}, which necessitates the evaluation of the Laplacian of $\mathcal{N}^{\BNAVEM}$. Since the dominant computational cost arises from this operation, a viable alternative would be to train ${\mathcal N}^{\BNAVEM}$ using a Variational Physics-Informed Neural Network (V-PINN) formulation \cite{KHARAZMI2021113547, Pintore2022, berrone2022solving}. V-PINNs require only first-order derivatives of the network, therefore reducing the computational burden. In this work, we focus on other aspects and defer the investigation of this extension to future work, as it primarily impacts training time.
\end{remark}

Concerning the memory usage, we observe that the high computational cost related to the computation of the Laplacian of the neural network requires the usage of a smaller architecture for B-NAVEM with respect to NAVEM or P-NAVEM. More precisely, we adopt an architecture consisting of $5$ layers and $50$ neurons per layer to train NAVEM and P-NAVEM for each class of polygon. The same architecture is also used to train B-NAVEM on the classes of convex and concave quadrilaterals ($1000$ elements each), and convex heptagons (about $300$ elements). To avoid out-of-memory problems, we have to use the smaller architecture consisting of $4$ layers and $40$ neurons per layer to train B-NAVEM using the same training datasets consisting of $1000$ elements for the classes of convex pentagons and hexagons. 

\begin{figure}[!ht]
\centering
\begin{subfigure}{0.49\textwidth}
    \includegraphics[width=\textwidth]{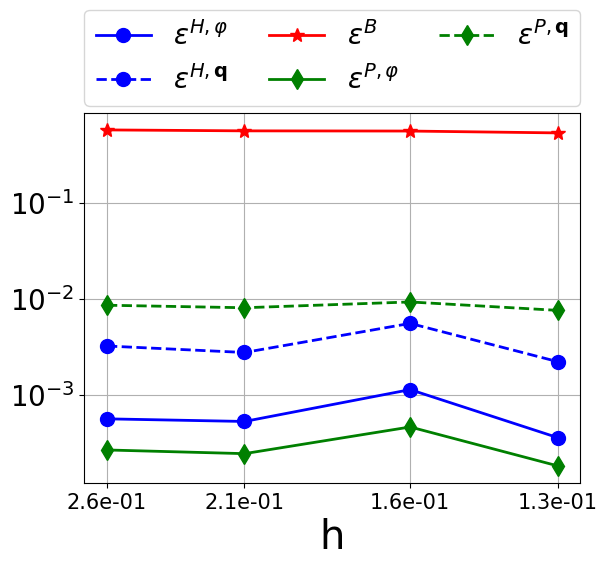}
    \caption{}
    \label{fig:loss_voro}
\end{subfigure}\hfill
\begin{subfigure}{0.49\textwidth}
    \includegraphics[width=\textwidth]{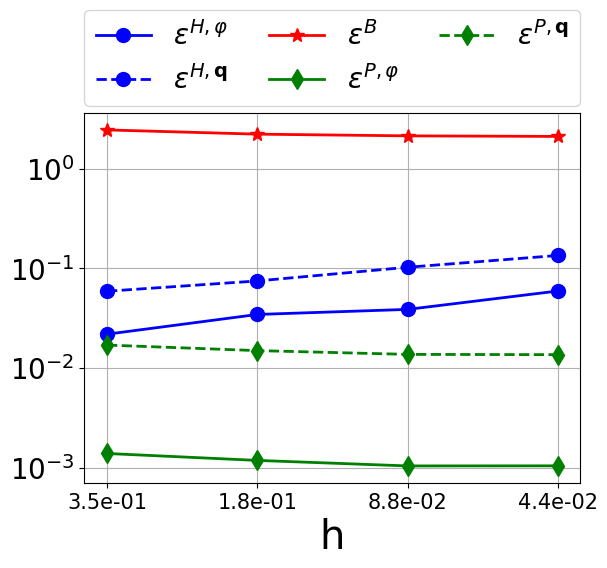}
    \caption{}
    \label{fig:loss_concave}
\end{subfigure}
\caption{Tests 1, 2, and 3: Test Loss values on each mesh refinement. Left: Voronoi. Right: Convex-Concave.}
\label{fig:compare_losses}
\end{figure}

Figure~\ref{fig:compare_losses} shows the behaviour of the test losses for each family of meshes and for each method. We note that, in this plot as well as in the following ones, the mesh size $h$ increases along the positive $x_1$-axis, i.e., the $x_1$-axis is inverted with respect the usual convention when it reports $h$ values. In particular, for the case of the P-NAVEM method, we report both the cost functionals \eqref{eq:napem:loss_basis} and \eqref{eq:napem:loss_gradients}, even if we remember that we train the P-NAVEM neural network just by minimizing the loss \eqref{eq:napem:loss_gradients} that involves gradients.
We can observe that since the training datasets we built are very representative of the test datasets, training and test performance are very similar to each other.

\begin{figure}[!ht]
\centering
\begin{subfigure}{0.49\textwidth}
    \includegraphics[width=\textwidth]{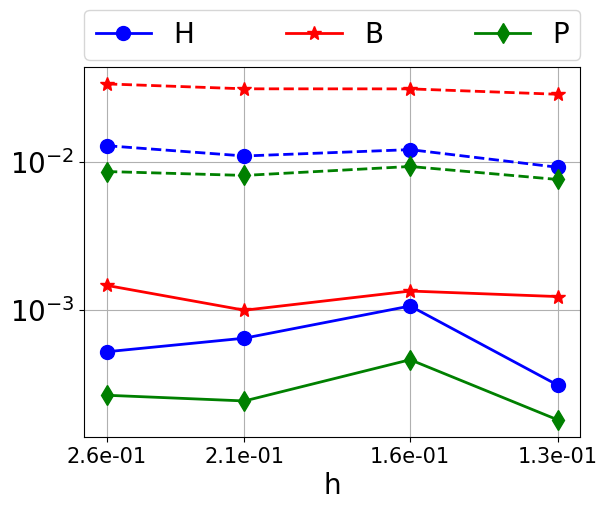}
    \caption{}
    \label{fig:pol_loss_voro}
\end{subfigure}\hfill
\begin{subfigure}{0.49\textwidth}
    \includegraphics[width=\textwidth]{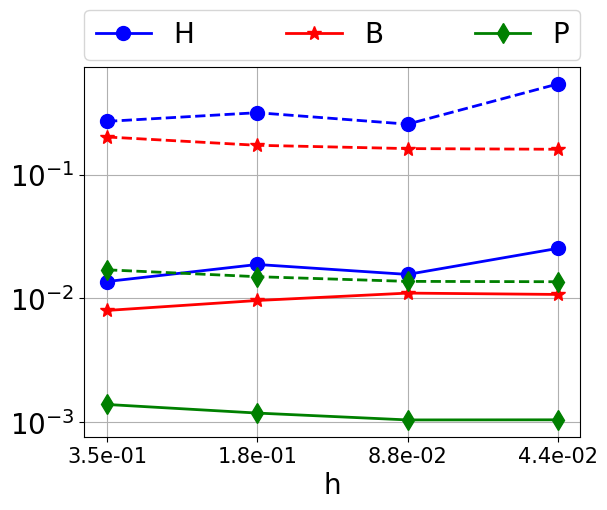}
    \caption{}
    \label{fig:pol_loss_concave}
\end{subfigure}
\caption{Tests 1, 2, and 3: Behaviour of the metrics \eqref{eq:napem:loss_basis} and \eqref{eq:napem:loss_gradients} computed on the test datasets for each method and for each test mesh refinement. Left: Voronoi. Right: Concave-convex.}
\label{fig:compare_pol_losses}
\end{figure}

Concerning the values of the various losses, we observe that the different nature of the loss functions associated with the different methods prevents a fair direct comparison.

To compare the accuracy of the different neural-based methods during the online test phase, it is therefore necessary to introduce a meaningful and common evaluation metric.

In particular, adopting the NAVEM or B-NAVEM loss functions as comparison metrics would be misleading, since these losses vanish identically for some networks by construction. More specifically, the NAVEM loss is zero for both the B-NAVEM and P-NAVEM networks, while the B-NAVEM loss is zero for the NAVEM network. As a result, these losses do not provide a meaningful metric for comparing the different methods.

For this reason, to enable a fair and informative comparison on the test datasets, we evaluate all neural-based methods using the loss functions \eqref{eq:napem:loss_basis} and \eqref{eq:napem:loss_gradients}, introduced in the P-NAVEM framework, as evaluation metrics. These losses can be consistently computed for all methods and represent the most relevant indicators of accuracy and convergence when numerically solving the Problem \eqref{eq:continuous_poisson}.

The results are reported in Figure \ref{fig:compare_pol_losses} for each test mesh family. Here, we can see that, in accordance with the fact that the P-NAVEM networks minimize an analogous cost function, their values are always smaller with respect to the other two methods. We refer to the next sections for further comments.

Finally, we want to highlight that the values of both the loss functions, shown in Figure \ref{fig:compare_losses}, and metrics \eqref{eq:napem:loss_basis}- \eqref{eq:napem:loss_gradients}, reported in Figure \ref{fig:compare_pol_losses}, are computed not only on polygons that are different from those contained in the training datasets, but also using a different set of points than the ones used in the training phase, in order to really test the generalization properties on these particular meshes. More precisely, these points are obtained using $N=13$ in the Algorithm \ref{alg:points_on_triangle}, instead of the value $N=10$ used in the training phase. 

\subsection{Test 2: A diffusion-advection-reaction problem}\label{sec:test2}
\begin{figure}[!ht]
    \centering
    \includegraphics[width=0.6\textwidth]{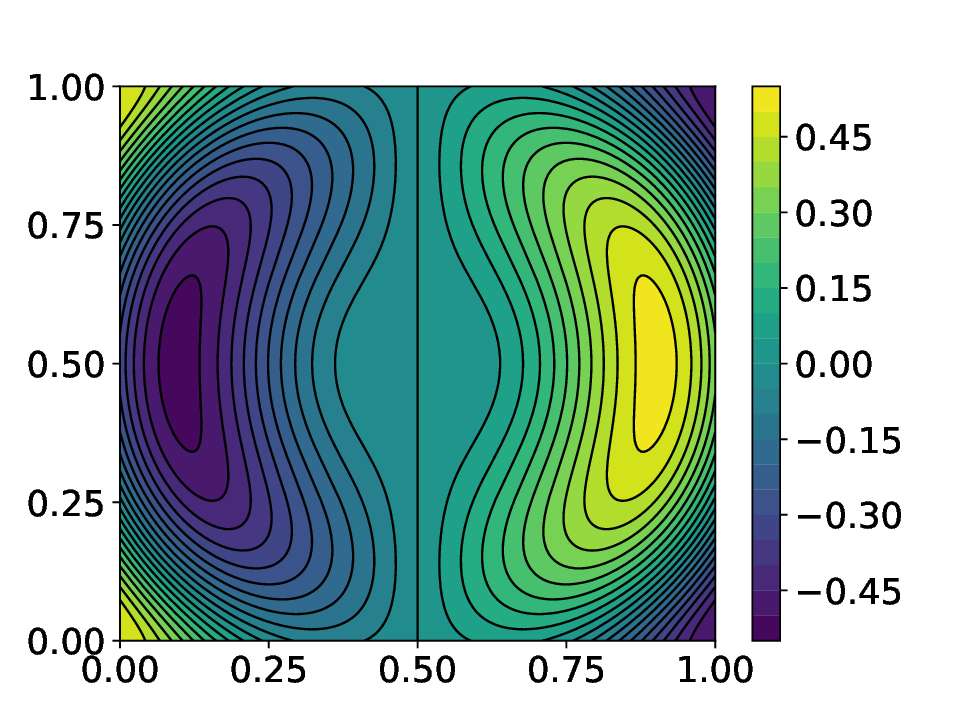}
    \caption{Test 2: Contour plot of the exact solution $u$ defined in \eqref{eq:test1_exact_solution}.}
    \label{fig:smb_solution}
\end{figure}

Let us consider the domain $\Omega=(0,1)^2$ with boundary $\partial\Omega$ and the following diffusion-advection-reaction problem:
\begin{equation}
    \begin{cases}
       - \nabla \cdot \left( \DD(\xx) \nabla u\right) + \bbeta(\xx) \cdot \nabla u + \gamma(\xx) u = f & \text{in } \Omega,\\
        u = g_D & \text{on } \partial\Omega,
    \end{cases}
    \label{eq:probl_convergence}
\end{equation} 
where 
\begin{linenomath}
\begin{equation*}
    \DD(\xx) = \begin{bmatrix}
        1 + x_2^2 & -x_1 x_2\\
        -x_1 x_2 & 1 + x_1^2
    \end{bmatrix}, \quad\quad \bbeta(\xx) = \begin{bmatrix}
     x_1 \\
     -x_2
    \end{bmatrix}, \quad\quad \gamma(\xx) = x_1 x_2,
\end{equation*}
\end{linenomath}
and where the forcing term $f$ and the Dirichlet boundary condition $g_D$ are chosen so that the exact solution, shown in Figure \ref{fig:smb_solution}, is
\begin{equation}
 u(\xx) = \sin\left[\pi^2 \left(\left(x_1-\dfrac12\right)^2 + \left(x_2-\dfrac12\right)^2\right)\right] \left(x_1 - \dfrac12\right) (1 + \sin(\pi x_1)  \sin(\pi x_2)).
     \label{eq:test1_exact_solution}
\end{equation} 

For comparison purposes, Problem \eqref{eq:probl_convergence} is solved with the standard VEM and NAVEM, and with the new methods B-NAVEM and P-NAVEM on the two mesh families $\{\mathcal{T}_{h,i}^Q\}_{i=1}^4$ and $\{\mathcal{T}_{h,i}^V\}_{i=1}^4$ described in the previous test. Moreover, we remark that, for neural-based methods, the neural networks employed in this simulation correspond to those described in the previous test case.

\begin{figure}[!ht]
    \centering
    \begin{subfigure}{0.49\textwidth}
        \includegraphics[width=\textwidth]{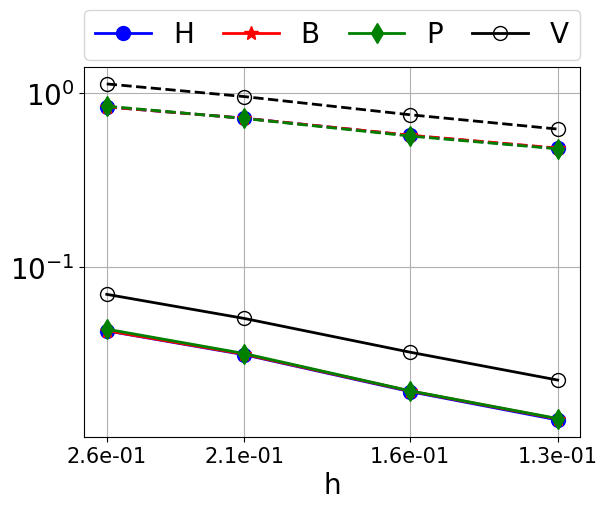}
        \caption{}
        \label{fig:error_voro}
    \end{subfigure}\hfill
    \begin{subfigure}{0.49\textwidth}
        \includegraphics[width=\textwidth]{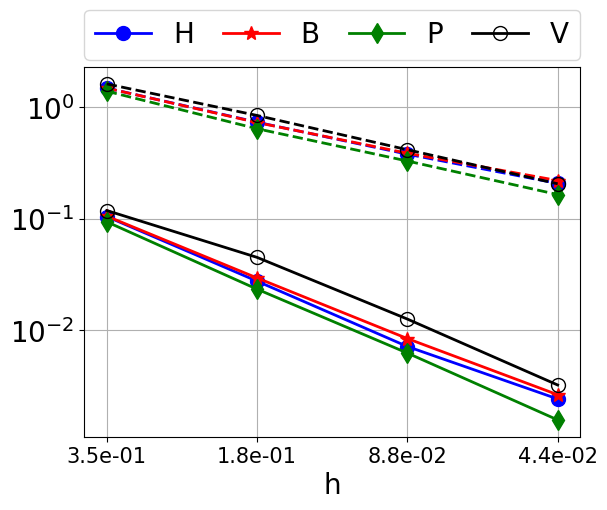}
        \caption{}
        \label{fig:error_concave}
    \end{subfigure}
    \caption{Test 2: Behaviour of the errors \eqref{eq:nn_errors} for the neural-based methods and of errors \eqref{eq:vem_errors} for the VEM method as $h$ decreases. Solid lines represent error $\mathrm{err}_0$, whereas dashed lines denote $\mathrm{err}_{\nabla}$ errors. Left: Voronoi. Rights: Convex-Concave.}
    \label{fig:compare_errors}
\end{figure}

Since the approximated solution obtained with NAVEM, B-NAVEM, and P-NAVEM is known in closed form, we measure the performance of these methods by looking at the behaviour of the following errors 
\begin{equation}
    \mathrm{err}_0^* = \sqrt{\sum_{E\in\Th} \norm[0,E]{u - u^{*}_h}^2},\quad \mathrm{err}_{\nabla}^* = \sqrt{\sum_{E\in\Th} \norm[0,E]{\nabla u - \nabla u^{*}_h}^2}, \quad *\in\{\NAVEM, \BNAVEM, \NAPEM\}.
    \label{eq:nn_errors}
\end{equation}
as the mesh size $h$ decreases. For the neural-based methods, the discrete solution and its gradient are given by
\begin{equation*}
    u_h^{\ast} = \sum_{i=1}^{\Ndof} u_i^{\ast} \varphi_{i}^{\ast}\quad \nabla u_h^{\ast} = \sum_{i=1}^{\Ndof} \qq_{i}^{\ast}\quad \forall \ast \in \{\NAVEM, \BNAVEM, \NAPEM\},
\end{equation*}
where the vector containing the degrees of freedom $\{u_i^{\ast}\}_{i=1}^{\Ndof}$ is the solution of the discrete Problem \eqref{eq:discrete_var_prob}.
On the other hand, since the VEM solution is not known in a closed form inside each element, we need to resort to polynomial projections of virtual functions to access their point-wise evaluation and compute the errors. Thus, we consider the usual definition of the VEM errors \cite{LBe16}:
\begin{equation}
    \mathrm{err}_0^{\VEM} = \sqrt{\sum_{E\in\Th} \norm[0,E]{u - \proj{0,E}{1} u_h^{\VEM}}^2},\quad \mathrm{err}_{\nabla}^{\VEM} = \sqrt{\sum_{E\in\Th} \norm[0,E]{\nabla u - \proj{\nabla,E}{1} \nabla u_h^{\VEM}}^2},
    \label{eq:vem_errors}
\end{equation}
where $\proj{0,E}{1}$ denotes the $\leb{2}{}$-projection onto $\Poly{1}{E}$, while $\proj{\nabla,E}{1}$ represents the elemental projection onto $\Poly{1}{E}$ in the energy norm.

The errors decays related to the two mesh families are shown in Figure~\ref{fig:compare_errors} for each method. It can be observed that the convergence rates of for neural-based methods are very close to the VEM ones in both cases. Moreover, we can notice that the error curves related to neural-based methods are downward shifted with respect the VEM ones for both the mesh families. We further note that all the neural-based methods perform very similarly on Voronoi meshes. On the other hand, P-NAVEM is the most accurate on meshes involving concave elements.

\begin{figure}[!ht]
\centering
\begin{subfigure}{0.49\textwidth}
    \includegraphics[width=\textwidth]{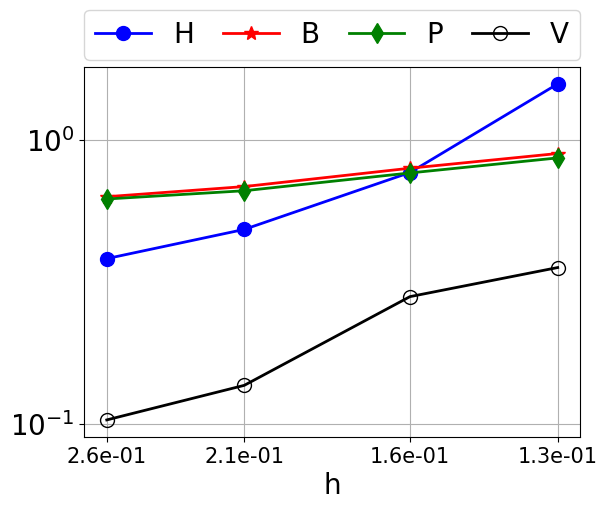}
    \caption{}
    \label{fig:assemble_time_voro}
\end{subfigure}\hfill
\begin{subfigure}{0.49\textwidth}
    \includegraphics[width=\textwidth]{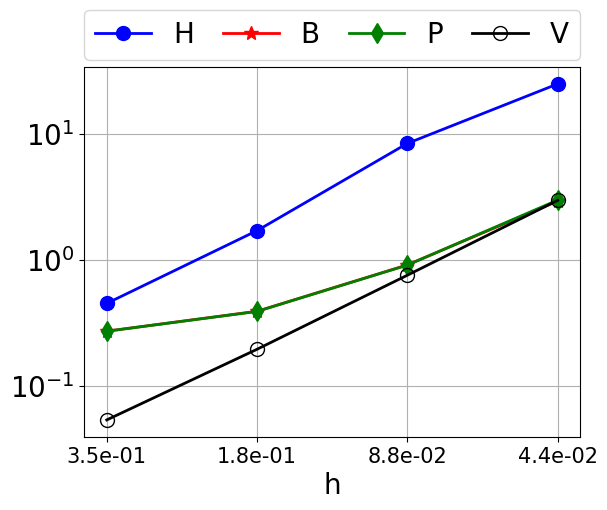}
    \caption{}
    \label{fig:assemble_time_concave}
\end{subfigure}
\caption{Test 2: Time required (in seconds) to assemble the global system matrix and solve the discrete problem associated with \eqref{eq:probl_convergence}. Left: Voronoi. Right: Convex-Concave.}
\label{fig:compare_assemble_time}
\end{figure}

To conclude this experiment, we compare the different methods by looking at the computational time $T$ (in seconds) required to solve Problem \eqref{eq:probl_convergence} for each mesh refinement with each method. In particular, the computed time $T$ represents an average time over $10$ different executions of a Python code that exploits TensorFlow \cite{tensorflow2015-whitepaper} for neural network operations, that runs on a Ubuntu 24.04 LTS 64-bit, 12th Gen Intel(R) Core(TM) i7-1255U CPU (4.7 GHz) and 16 GB RAM memory.

More specifically, for each neural-based method, this time accounts for 
\begin{itemize}
    \item the time needed to evaluate auxiliary functions at quadrature points, like the harmonic functions in $\approxspace[j,E]$ and their derivatives for NAVEM, the ADF functions \eqref{eq:adf_boundary}, the transfinite interpolant \eqref{eq:transfite_interpolant}, and their derivatives for B-NAVEM and P-NAVEM;
    \item the time needed to encode the neural network input and to produce the output;
    \item the time to assemble and solve the associated discrete system.
\end{itemize}
For the standard VEM, instead, the computational time $T$ accounts for the time needed to compute projection and stability operators as well as time needed to assemble and solve the global discrete system.

The behaviour of the computational time $T$ for each method and mesh family is shown in Figure~\ref{fig:compare_assemble_time}. From these log–log plots, we observe that neural-based methods generally require more computational time than the standard VEM. However, such a difference decreases when employing the B-NAVEM or the P-NAVEM method for finer meshes.

We remark that, to mitigate the cost associated with neural network evaluations, we minimize the number of neural network evaluations by aggregating the inputs associated with multiple basis functions, thus exploiting TensorFlow vectorization capabilities and reducing the overall computational burden. Nevertheless, the intrinsic cost per network call has a stronger impact on small datasets, corresponding to coarse meshes, and on mesh families characterized by a high variability in the number of polygon vertices (polygon classes), as already observed in \cite{PintoreTeora2025}.

No significant difference is observed between the computational costs of the B-NAVEM and P-NAVEM methods. Indeed, both approaches require the same number of neural network evaluations and the same computational effort during the test phase to evaluate the auxiliary functions needed to enforce boundary conditions. Moreover, we observe that this overhead with respect to VEM becomes negligible as $h$ decreases for P-NAVEM and B-NAVEM, since these auxiliary functions can be computed very efficiently using TensorFlow and NumPy vectorization \cite{numpy}.

On the other hand, the computation of the harmonic functions required to evaluate the NAVEM basis functions is performed element-wise. About this, we observe that we always adopt a value for the harmonic polynomial degree $\ell^{\NN} = 20$ for all the NAVEM neural networks we train, without fine-tuning this value, further increasing the cost associated with harmonic function evaluations. As a consequence, the associated computational cost increases linearly with the number of elements.

In conclusion, as the mesh size $h$ decreases, even for a simple advection–diffusion–reaction problem, the neural-based methods become competitive with the standard VEM in terms of computational effort, while providing higher accuracy, up to the approximation capability of the neural network. 

\subsection{Test 3: A nonlinear diffusion problem}

\begin{figure}[!ht]
    \centering
    \begin{subfigure}{0.32\textwidth}
        \includegraphics[width=\textwidth]{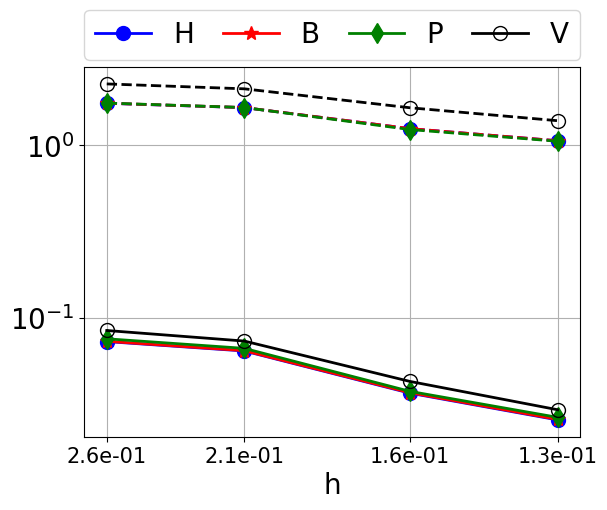}
        \caption{}
        \label{fig:error_voro_nl_ring_1}
    \end{subfigure}\hfill
    \begin{subfigure}{0.32\textwidth}
        \includegraphics[width=\textwidth]{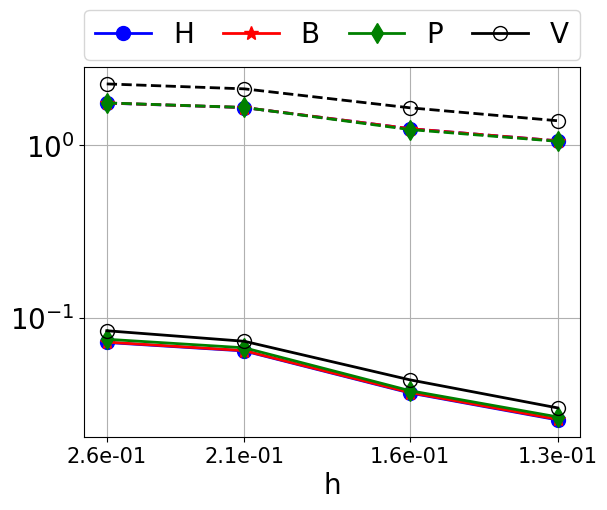}
        \caption{}
        \label{fig:error_voro_nl_ring_05}
    \end{subfigure}\hfill
    \begin{subfigure}{0.32\textwidth}
        \includegraphics[width=\textwidth]{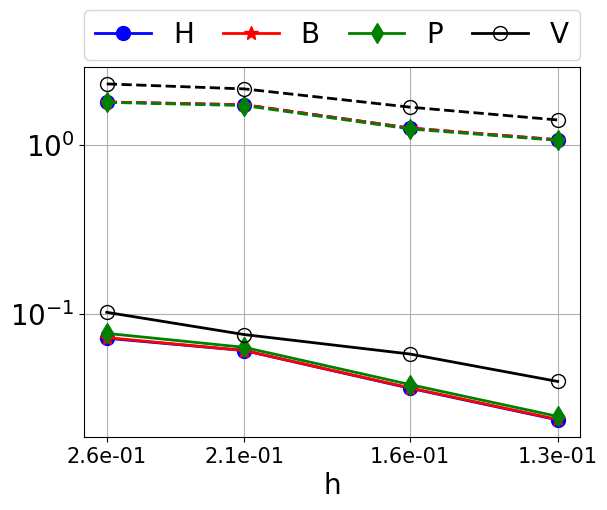}
        \caption{}
        \label{fig:error_voro_nl_ring_01}
    \end{subfigure}
    \begin{subfigure}{0.32\textwidth}
        \includegraphics[width=\textwidth]{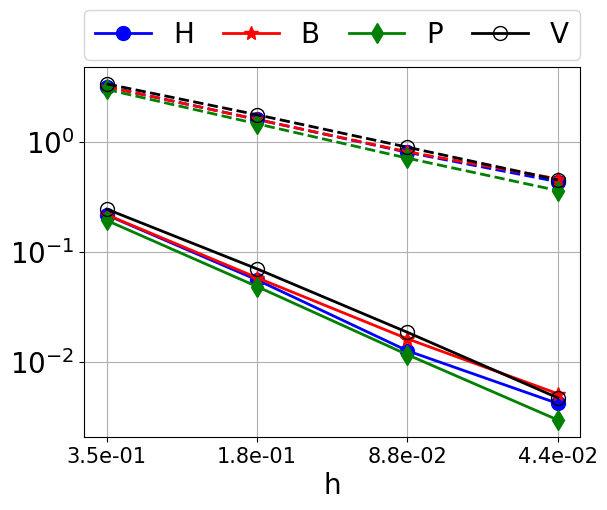}
        \caption{}
        \label{fig:error_concave_nl_ring_1}
    \end{subfigure}
    \begin{subfigure}{0.32\textwidth}
        \includegraphics[width=\textwidth]{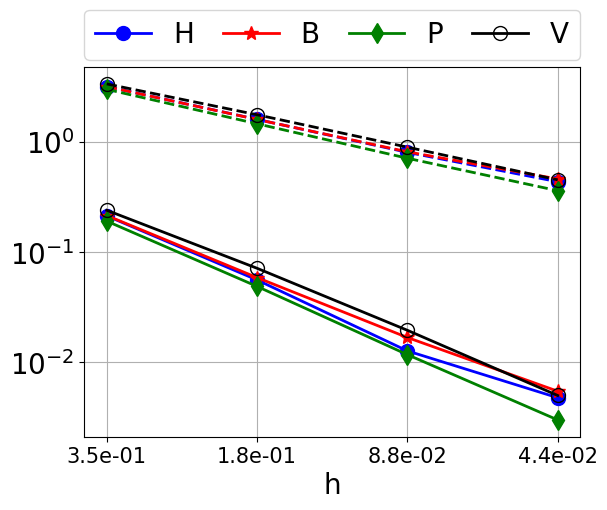}
        \caption{}
        \label{fig:error_concave_nl_ring_05}
    \end{subfigure}
    \begin{subfigure}{0.33\textwidth}
        \includegraphics[width=\textwidth]{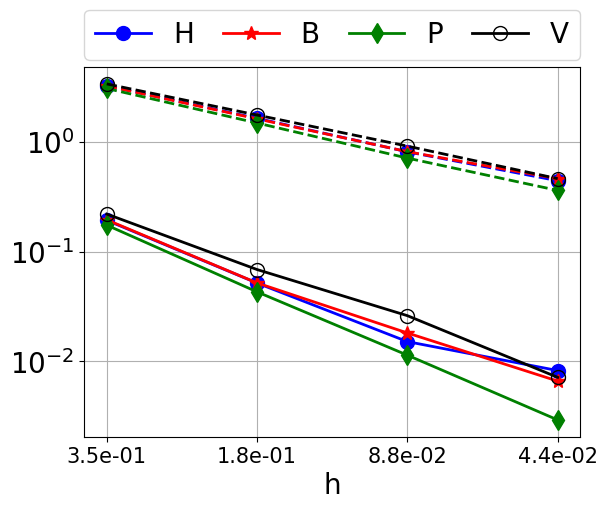}
        \caption{}
        \label{fig:error_concave_nl_ring_01}
    \end{subfigure}
    \caption{Test 3: Behaviour of the errors \eqref{eq:nn_errors} for the neural-based methods and of errors \eqref{eq:vem_errors} for the VEM method as $h$ decreases. Solid lines represent error $\mathrm{err}_0$, whereas dashed lines denote $\mathrm{err}_{\nabla}$ errors. Each column refers to a value of $\lambda$, namely $1.0$, $0.5$, and $0.1$ from left to right. First row: Voronoi. Second row: Convex-Concave.}
    \label{fig:error_nl_ring}
\end{figure}

\begin{figure}[!ht]
    \centering
    \begin{subfigure}{0.32\textwidth}
        \includegraphics[width=\textwidth]{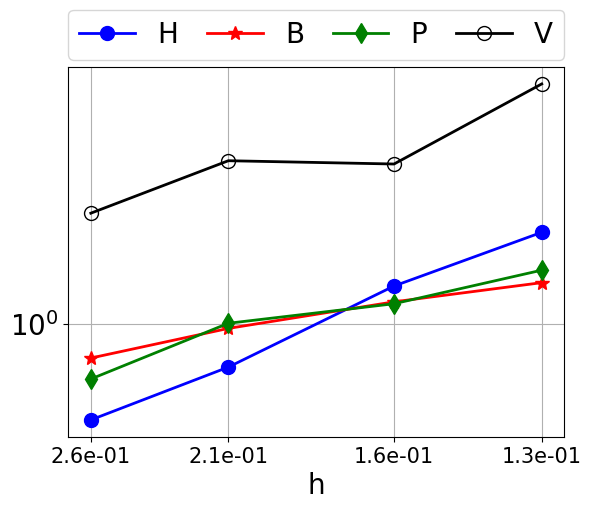}
        \caption{}
        \label{fig:assemble_time_voro_nl_ring_1}
    \end{subfigure}\hfill
    \begin{subfigure}{0.32\textwidth}
        \includegraphics[width=\textwidth]{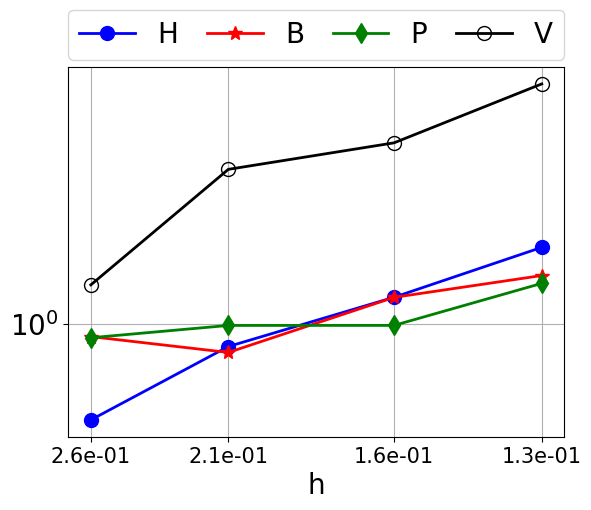}
        \caption{}
        \label{fig:assemble_time_voro_nl_ring_05}
    \end{subfigure}\hfill
    \begin{subfigure}{0.32\textwidth}
        \includegraphics[width=\textwidth]{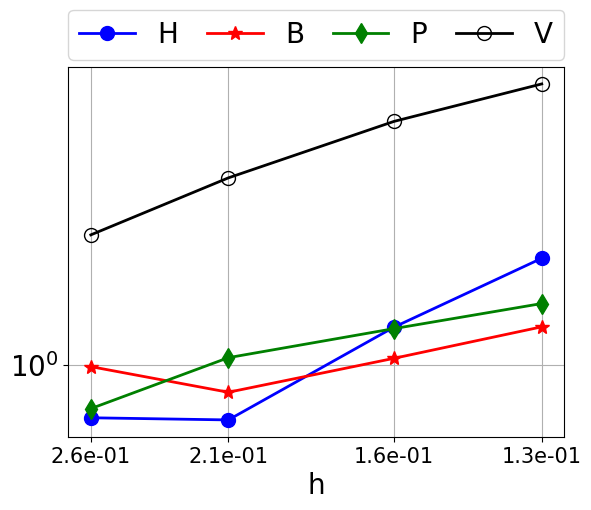}
        \caption{}
        \label{fig:assemble_time_voro_nl_ring_01}
    \end{subfigure}
    \begin{subfigure}{0.32\textwidth}
        \includegraphics[width=\textwidth]{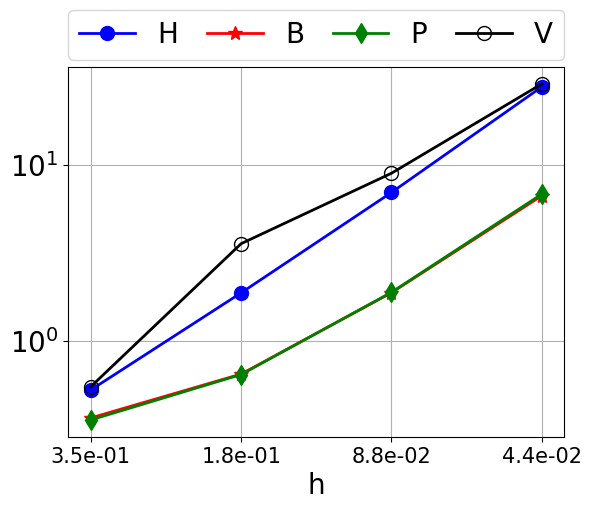}
        \caption{}
        \label{fig:assemble_time_concave_nl_ring_1}
    \end{subfigure}
    \begin{subfigure}{0.32\textwidth}
        \includegraphics[width=\textwidth]{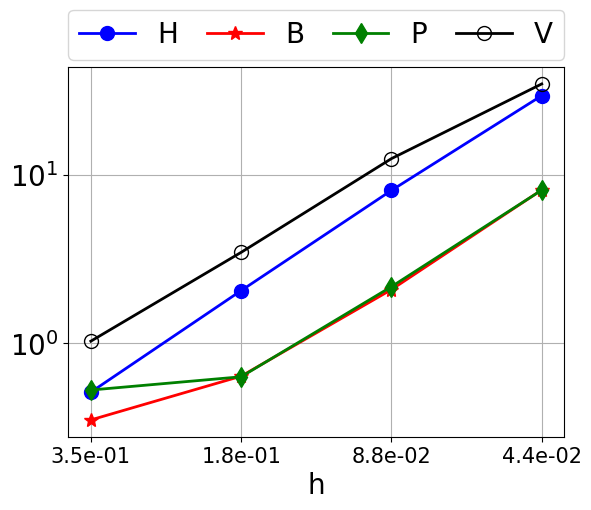}
        \caption{}
        \label{fig:assemble_time_concave_nl_ring_05}
    \end{subfigure}
    \begin{subfigure}{0.32\textwidth}
        \includegraphics[width=\textwidth]{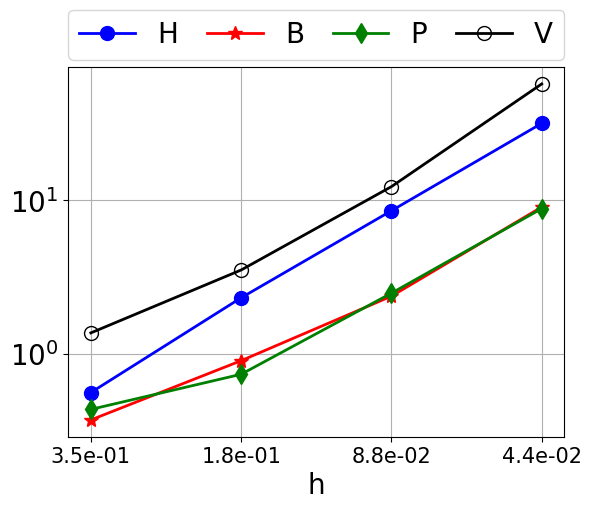}
        \caption{}
        \label{fig:assemble_time_concave_nl_ring_01}
    \end{subfigure}
    \caption{Test 3: Time required (in seconds) to assemble the global system matrix and solve the discrete problem associated with \eqref{eq:nn_diffusion_problem}. Each column refers to a value of $\lambda$, namely $1.0$, $0.5$, and $0.1$ from left to right. First row: Voronoi. Second row: Convex-Concave.}
    \label{fig:compare_assemble_time_nl_ring}
\end{figure}

To show the benefits of avoiding projection and stability operators that are typical of the standard virtual element methods, in this section, we replicate Test Problem 3 in \cite{PintoreTeora2025} on the new families of meshes $\{\mathcal{T}_{h,i}^V\}_{i=1}^4$ and $\{\mathcal{T}_{h,i}^Q\}_{i=1}^4$. More precisely, we set $\Omega=(0,1)^2$ and consider the following nonlinear diffusion problem
\begin{equation}
    \begin{cases}
       - \nabla \cdot \left( \DD(u, \lambda) \nabla u\right)  = f & \text{in } \Omega,\\
        u = g_D & \text{on } \partial\Omega,
    \end{cases}
    \label{eq:nn_diffusion_problem}
\end{equation} 
where $\lambda \in \{1.0, 0.5, 0.1\}$ is a problem parameter, and the diffusion coefficient is defined as
\begin{equation*}
    \DD(u, \lambda) = \frac{1}{\lambda + u^2}.
\end{equation*}

For computing the errors \eqref{eq:nn_errors} and \eqref{eq:vem_errors}, we set the Dirichlet boundary condition and the forcing term in such a way that the exact solution is
\begin{equation*}
    u(\xx) = \frac{1}{8}\Big[\sin\left(3\pi((x_1 - 0.5)^2 + (x_2 -0.5)^2)\right)\Big]^3.
\end{equation*}

As in \cite{PintoreTeora2025}, to solve the nonlinear Problem \eqref{eq:nn_diffusion_problem}, we adopt the Newton method, and, by setting as the initial guess the all zero vector, we adopt the following stop criteria:
\begin{equation*}
    \norm{\bm{r}_h^{m,\ast}} \leq \alpha_{r,f} \norm{\bm{f}_h^{0,\ast}} + \alpha_{a,f} \quad \text{and} \quad \norm{\bm{\delta}_h^{m,\ast}} \leq \alpha_{a,\delta} \quad \forall \ast \in \{\VEM,\NAVEM,\BNAVEM,\NAPEM\},
\end{equation*}
where $\bm{r}_h^{m,\ast}$ is the residual associated with the method at the nonlinear iteration $m$, $\bm{f}_h^{0,\ast}$ is the initial residual, $\bm{\delta}_h^{m,\ast}$ is the difference between the discrete solutions at the iteration $m + 1$ and $m$, and $\alpha$ values are tolerance chosen by the user.

In this test, to compare the different methods, we compute the errors defined in \eqref{eq:nn_errors} and \eqref{eq:vem_errors} to assess their accuracy. In addition, to highlight the advantages of avoiding projection and stability operators when solving nonlinear problems, we introduce the following performance indicators. For each neural-based method $\ast \in \{\NAVEM,\BNAVEM,\NAPEM\}$, we compute:
\begin{enumerate}[label=\textbf{S.\arabic*}]
    \item \label{stat:time1}$r^{\ast}(T)$: the ratio between the time $T$ (in seconds) required by VEM and by method $\ast$ to solve the nonlinear Problem \eqref{eq:nn_diffusion_problem};
    \item $r^{\ast}(m)$: the ratio between the number of nonlinear iteration $m$ required by VEM and by method $\ast$;
    \item \label{stat:time3}$r^{\ast}(\rm{ATI})$: the ratio between the Average Time per Iteration (ATI) in VEM and in approach $\ast$.
\end{enumerate}

In Figures~\ref{fig:error_nl_ring} and~\ref{fig:compare_assemble_time_nl_ring}, we report the error behaviour and the total computational time $T$, computed as described in Section~\ref{sec:test2}, respectively, as the mesh size $h$ decreases, for each method and mesh family.

As already observed in Test~2 (Section~\ref{sec:test2}), avoiding projection and stability operators leads to smaller error constants compared to the standard VEM. Moreover, the P-NAVEM method, characterized by the smaller polynomial losses \eqref{eq:napem:loss_basis} and \eqref{eq:napem:loss_gradients}, is confirmed to be the most accurate among the neural-based methods for sufficiently small values of $h$, especially in the presence of concave elements. Furthermore, we observe that as the problem becomes more nonlinear (i.e., as $\lambda$ decreases), the accuracy advantage of neural-based methods over VEM becomes more pronounced.

\begin{table}[!ht]
\centering
\caption{Test 3: Statistics about computational time required by each method to solve the nonlinear Problem \eqref{eq:nn_diffusion_problem} via Newton method for each mesh refinement $i$. Voronoi family of meshes. For each row, the values for the average time $T$ and the average time per iteration $\rm{ATI}$ related to the method with the best performance are highlighted in bold.}
\label{tab:time_voronoi}
\resizebox{0.95\textwidth}{!}{
\begin{tabular}{cc|ccc|ccc|ccc}
$\lambda$ & $i$ & $r^{\NAVEM}(T)$ & $r^{\NAVEM}(m)$ & $r^{\NAVEM}(\rm{ATI})$ & $r^{\BNAVEM}(T)$ & $r^{\BNAVEM}(m)$ & $r^{\BNAVEM}(\rm{ATI})$ & $r^{\NAPEM}(T)$ & $r^{\NAPEM}(m)$ & $r^{\NAPEM}(\rm{ATI})$ \\ [0.5mm] \hline\addlinespace[0.5mm]
1.0       & 0   & \textbf{3.15}            & 1.67            & \textbf{1.89}                   & 2.23             & 1.67             & 1.34                    & 2.51            & 1.67            & 1.51                   \\
0.5       & 0   & \textbf{2.56}            & 2.00            & \textbf{1.28}                   & 1.43             & 2.00             & 0.72                    & 1.45            & 2.00            & 0.72                   \\
0.1       & 0   & \textbf{3.30}            & 1.75            & \textbf{1.88}                   & 2.36             & 1.75             & 1.35                    & 3.11            & 1.75            & 1.77                   \\ [1mm]
1.0       & 1   & \textbf{3.14}            & 1.67            & \textbf{1.89}                   & 2.54             & 1.67             & 1.52                    & 2.46            & 1.67            & 1.48                   \\
0.5       & 1   & 3.44            & 1.83            & 1.87                   & \textbf{3.58}             & 1.83             & \textbf{1.95}                    & 2.97            & 1.83            & 1.62                   \\
0.1       & 1   & \textbf{4.85}            & 1.75            & \textbf{2.77}                   & 4.05             & 1.75             & 2.31                    & 3.23            & 1.75            & 1.85                   \\ [1mm]
1.0       & 2   & 1.97            & 1.67            & 1.18                   & \textbf{2.15}             & 1.67             & \textbf{1.29}                    & 2.17            & 1.67            & 1.30                   \\
0.5       & 2   & 2.93            & 1.83            & 1.60                    & 2.94            & 1.83            & 1.60                    & \textbf{3.57}            & 1.83            & \textbf{1.95}                   \\
0.1       & 2   & 3.83            & 1.50            & 2.55                   & \textbf{4.69}             & 1.50             & \textbf{3.13}                    & 3.86            & 1.50            & 2.58                   \\ [1mm]
1.0       & 3   & 2.28            & 1.50            & 1.52                   & \textbf{3.01}             & 1.50             & \textbf{2.01}                    & 2.81            & 1.50            & 1.87                   \\
0.5       & 3   & 3.12            & 1.67            & 1.87                   & 3.80             & 1.67             & 2.28                    & \textbf{4.01}            & 1.67            & \textbf{2.41}                   \\
0.1       & 3   & 3.12            & 1.50            & 2.08                   & \textbf{4.88}             & 1.50             & \textbf{3.25}                    & 4.19            & 1.50            & 2.79                   \\ \hline
\end{tabular}}
\end{table}

\begin{table}[!ht]
\centering
\caption{Test 3: Statistics about computational time required by each method to solve the nonlinear Problem \eqref{eq:nn_diffusion_problem} for each mesh refinement $i$. Convex-Concave family of meshes. For each row, the values for the average time $T$ and the average time per iteration $\rm{ATI}$ related to the method with the best performance are highlighted in bold.}
\label{tab:time_quad}
\resizebox{0.95\textwidth}{!}{
\begin{tabular}{cc|ccc|ccc|ccc}
$\lambda$ & $i$ & $r^{\NAVEM}(T)$ & $r^{\NAVEM}(m)$ & $r^{\NAVEM}(\rm{ATI})$ & $r^{\BNAVEM}(T)$ & $r^{\BNAVEM}(m)$ & $r^{\BNAVEM}(\rm{ATI})$ & $r^{\NAPEM}(T)$ & $r^{\NAPEM}(m)$ & $r^{\NAPEM}(\rm{ATI})$ \\ [0.5mm] \hline\addlinespace[0.5mm]
1.0       & 0   & 1.04            & 1.50            & 0.70                   & 1.51             & 1.50             & 1.00                    & \textbf{1.55}            & 1.50            & \textbf{1.03}                   \\
0.5       & 0   & 2.01            & 1.83            & 1.10                   & \textbf{2.94}             & 1.83             & \textbf{1.60}                    & 1.95            & 1.83            & 1.06                   \\
0.1       & 0   & 2.45            & 2.57            & 0.95                   & \textbf{3.70}             & 2.57             & \textbf{1.44}                    & 3.15            & 2.57            & 1.23                   \\ [1mm]
1.0       & 1   & 1.90            & 1.67            & 1.14                   & 5.52             & 1.67             & 3.31                    & \textbf{5.55}            & 1.67            & \textbf{3.33}                   \\
0.5       & 1   & 1.69            & 1.83            & 0.92                   & 5.44             & 1.83             & 2.97                    &\textbf{ 5.48}            & 1.83            & \textbf{2.99}                   \\
0.1       & 1   & 1.52            & 1.75            & 0.87                   & 3.90             & 1.75             & 2.23                    &\textbf{4.79}            & 1.75            & \textbf{2.74}                   \\ [1mm]
1.0       & 2   & 1.29            & 1.50            & 0.86                   & \textbf{4.78}             & 1.50             & \textbf{3.19}                    & 4.76            & 1.50            & 3.17                   \\ 
0.5       & 2   & 1.54            & 1.43            & 1.08                   & \textbf{6.00}             & 1.67             & \textbf{3.60}                    & 5.76            & 1.67            & 3.46                   \\
0.1       & 2   & 1.44            & 1.50            & 0.96                   & \textbf{5.19}             & 1.50             & \textbf{3.46}                    & 4.96            & 1.50            & 3.30                   \\ [1mm]
1.0       & 3   & 1.04            & 1.33            & 0.78                   & \textbf{4.33}             & 1.33             & \textbf{3.25}                    & 4.24            & 1.33            & 3.18                   \\
0.5       & 3   & 1.17            & 1.29            & 0.91                   & 4.26             & 1.29             & 3.32                    & \textbf{4.27}            & 1.29            & \textbf{3.32}                   \\
0.1       & 3   & 1.81            & 1.25            & 1.45                   & 6.36             & 1.25             & 5.08                    & \textbf{6.51}            & 1.25            & \textbf{5.21}                   \\ \hline
\end{tabular}}
\end{table}

Figure~\ref{fig:compare_assemble_time_nl_ring} shows that the computational-time behaviour of neural-based methods is consistent with the one observed in the previous test case: for the coarser mesh and the Voronoi family, NAVEM is faster than B-NAVEM or P-NAVEM, but it is slower when simulating over finer meshes or over meshes belonging to Convex-Concave family. However, in this nonlinear setting, neural-based methods require significantly less computational time than VEM in total to solve Problem~\eqref{eq:nn_diffusion_problem}. This behaviour can be explained by two main factors. First, the overhead associated with neural-network utilities is amortized over all nonlinear iterations, thus improving the overall efficiency of neural-based methods. Second, these methods typically require fewer nonlinear iterations than VEM to achieve the desired accuracy, resulting in a smaller number of total calls to the direct solver, which is used to solve the linearized discrete problem at each iteration.

These observations are summarized by the statistics~\ref{stat:time1}–\ref{stat:time3} reported in Tables~\ref{tab:time_voronoi} and~\ref{tab:time_quad} for the Voronoi and Convex–Concave mesh families, respectively. From these tables, we highlight that neural-based methods employ the same number of nonlinear iterations to reach the desiderate accuracy, i.e. $r^{\NAVEM}(m) = r^{\BNAVEM}(m) = r^{\NAPEM}(m)$. We further observe that, for all mesh families, the best performance of neural-based methods over VEM are obtained for the smallest value of $\lambda=0.1$. This suggests that approaches that avoid stabilization and projection operators become increasingly efficient as the underlying PDE exhibits stronger nonlinearities.

\section{Conclusion}\label{sec:conclusion}

In this manuscript, we propose two polygonal discretization methods, called B-NAVEM and P-NAVEM, as alternatives to the standard NAVEM, in which local basis functions are constructed using pre-trained neural networks. In contrast to the standard NAVEM approach, the basis functions generated by the proposed methods are exactly continuous across adjacent elements, and this continuity is enforced by construction.

The B-NAVEM basis functions are defined through a Physics-Informed Neural Network that minimizes the residual of the elemental Laplace problems characterizing the virtual element basis functions. As in the classical VEM and contrary to NAVEM, B-NAVEM functions are linear polynomials on the boundary. However, their Laplacian is only approximately zero, rather than identically zero as in the VEM or NAVEM setting.

The P-NAVEM approach, on the other hand, does not aim at approximating the VEM space itself. Instead, it directly constructs a local approximation space that is exactly $\con{0}{}$-conformed and accurately reproduces polynomials, a property that is essential to guarantee optimal polynomial convergence rates.

A series of numerical experiments is presented to compare the performance of the different methods on both linear and nonlinear benchmark problems. The results indicate that the P-NAVEM method provides the best trade-off between accuracy and computational cost, both in the training stage and in the testing phase.

Finally, using the same cost function employed in P-NAVEM methods, in \cite{TeoraNeva2025} the authors devise a standard method, termed the Zipped Finite Element Method, that locally solves the same optimization problem to build a space of higher-order basis functions known in closed-form on star-shaped polygons. In this case, since the optimization problem is solved element-wise, the local Z-FEM space exactly contains polynomials to preserve the optimal order of convergence for any method order $k \geq 1$. Such a new method does not employ neural networks, thus its accuracy is not limited by the neural network accuracy, but we need to solve a local optimization problem for each element in the tessellation.

\section*{Acknowledgements}

We thank Professor Claudio Canuto from Politecnico di Torino for the valuable suggestions.

\bigskip

The author S.B. kindly acknowledges partial financial support provided by European Union through project Next Generation EU, M4C2, PRIN 2022 PNRR project P2022BH5CB\_001 ``Polyhedral Galerkin methods for engineering applications to improve disaster risk forecast and management: stabilization-free operator-preserving methods and optimal stabilization methods'', and by PNRR M4C2 project of CN00000013 National Centre for HPC, Big Data and Quantum Computing (HPC) (CUP: E13C22000990001). The author M.P. kindly acknowledges financial support provided by PEPR/IA (\url{https://www.pepr-ia.fr/}). The author G.T. kindly acknowledges the financial support provided by project NODES which has received funding from the MUR-M4C2 1.5 of PNRR funded by the European Union - NextGenerationEU (Grant agreement no. ECS00000036) and by the European Union through PRIN project 20227K44ME ``Full and Reduced order modelling of coupled systems: focus on non-matching methods and automatic learning (FaReX)'' (CUP: E53D23005510006).

\bibliographystyle{IEEEtranDOI}
\bibliography{biblio.bib}

\end{document}